\documentclass[11pt]{article}
\pdfoutput=1

\usepackage{bm}
\usepackage[in]{fullpage}  

\usepackage{url}
\usepackage{algorithm, algorithmicx,algpseudocode}
\usepackage{multirow} 
\usepackage{hhline}  
\usepackage{amsmath}
\usepackage{graphicx,epstopdf} 
\usepackage{amssymb}
\usepackage{amsthm}
\usepackage{color}
\usepackage{subfigure}
\usepackage{authblk}

\allowdisplaybreaks
\newcommand\numberthis{\addtocounter{equation}{1}\tag{\theequation}}

\numberwithin{equation}{section}
\DeclareMathOperator{\rank}{\mathrm{rank}}

\def\lb{\left(}
\def\rb{\right)}
\def\lcb{\left\{}
\def\rcb{\right\}}

\def\lsb{\left[}
\def\rsb{\right]}

\def\C{\mathbb{C}}

\def\P{\mathcal{P}}
\def\T{\mathcal{T}}
\def\H{\mathcal{H}}

\def\R{\mathbb{R}}

\def\M{\mathcal{M}}

\def\I{\mathcal{I}}
\def\BZ{\bm{Z}}

\def\BL{\bm{L}}
\def\BU{\bm{U}}
\def\BV{\bm{V}}
\def\BB{\bm{B}}
\def\BC{\bm{C}}
\def\BW{\bm{W}}
\def\BD{\bm{D}}

\def\BM{\bm{M}}
\def\BQ{\bm{Q}}
\def\BH{\bm{H}}
\def\BX{\bm{X}}
\def\BY{\bm{Y}}
\def\BA{\bm{A}}
\def\BB{\bm{B}}
\def\BC{\bm{C}}
\def\BD{\bm{D}}
\def\BG{\bm{G}}
\def\BR{\bm{R}}

\def\BOMG{\bm{\Omega}}
\def\by{\bm{X}}
\def\BS{\bm{\Sigma}}
\def\BI{\bm{I}}

\def\bx{\bm{x}}
\def\by{\bm{y}}
\def\bz{\bm{z}}
\def\be{\bm{e}}
\def\bw{\bm{w}}

\begin{document}
\title{Fast Cadzow's Algorithm and a Gradient Variant\footnotetext{Corresponding author: Ke Wei (kewei@fudan.edu.cn).}}
\author[1]{Haifeng Wang}
\author[2]{Jian-Feng Cai}
\author[3]{Tianming Wang}
\author[1]{Ke Wei}
\affil[1]{School of Data Science, Fudan University, Shanghai, China.\vspace{.15cm}}
\affil[2]{Department of Mathematics, Hong Kong University of Science and Technology, Clear Water Bay, Kowloon, Hong Kong SAR, China.\vspace{.15cm}}
\affil[3]{Oden Institute of Computational Engineering and Sciences, The University of Texas at Austin, Austin, Texas, USA.}
\date{}
\maketitle
\begin{abstract}
The Cadzow's algorithm is a signal denoising and recovery method which was designed for signals corresponding to low rank Hankel matrices.   In this paper we first introduce a Fast Cadzow's algorithm which is developed by incorporating a novel subspace projection to reduce the high computational cost of the SVD in the Cadzow's algorithm. Then a Gradient method and a Fast Gradient method are proposed to address the non-decreasing MSE issue when applying the Cadzow's or Fast Cadzow's algorithm for signal denoising. Extensive empirical performance comparisons demonstrate that the proposed algorithms can complete the denoising and recovery tasks more efficiently and effectively.
\end{abstract}
\section{Introduction}\label{sec:introduction}
Cadzow's algorithm \cite{cadzow1988signal} was proposed by J. A. Cadzow  in 1988. The algorithm was initially developed for signal denoising, but has already been extended to many other  applications, including time series forecasting \cite{golyandina2005ssa}, the filtering of digital terrain models  \cite{golyandina2007filtering}, and seismic data denoising and reconstruction \cite{sacchi2009fx,oropeza2009multifrequency,trickett2008f,trickett2009prestack,oropeza2011simultaneous,CGW19}. For simplicity, we consider the 1D signal denoising problem. Here the task is to estimate a target signal $\bx\in\C^N$ from the corrupted measurements
\begin{align*}
\by = \bx+\be,\numberthis\label{eq:denoising}
\end{align*}
where $\be\in\C^N$ denotes additive noise. It is not hard to see that, without any additional assumptions, one cannot expect a universally better estimator than the one simply given by the noisy vector $\by$. Thus, effective denoising methods
typically rely on certain inherent simple structures hidden in  the target signal. For example, the method of wavelet denoising assumes that the representation of $\bx$ under certain wavelet transform has many entries that are close to zero \cite{donoho1995denoising,donohoJohn1998wavelet}. By contrast,  the method of Cadzow's denoising  is based on the low rank property of the Hankel  matrix associated with the target signal.

Let $\bz=[z_0,\cdots,z_{N-1}]^\top$ be a complex-valued vector of length $N$. Let $\H$ be a linear operator which maps $\bz$ into an $L\times K$   $(L+K=N+1)$ matrix $\H\bz$  whose $(i,j)$-th entry is equal to $z_{i+j}$, i.e.,
\begin{align*}
\H\bz = \begin{bmatrix}
z_0 & z_1 & z_2 & \cdots & z_{K-1}\\
z_1 & z_2 & z_3 & \cdots & z_{K}\\
z_2 & z_3 & z_4 & \cdots & z_{K+1}\\
\vdots & \vdots & \vdots & \vdots & \vdots\\
z_{L-1} & z_{L} & z_{L+1} & \cdots & z_{N-1}
\end{bmatrix},\numberthis\label{eq:hankel}
\end{align*}
where the numbering of vector and matrix entries starts from $0$ rather than the more standard $1$.  Matrices of the form \eqref{eq:hankel} are known as Hankel matrices in which each  skew-diagonal is  constant. Here we will refer to the Hankel matrix $\H\bz$ as the Hankelization of $\bz$. Furthermore, the Moore-Penrose
pseudoinverse of $\H$, denoted $\H^\dagger$, is a linear map from $L\times K$ matrices to  vectors of length $N$. 
Given a matrix $\BZ\in\C^{L\times K}$, the vector $\H^\dagger\BZ\in\C^N$ is obtained by averaging each skew-diagonal of $\BZ$. More precisely, letting $w_a$ $(a=0,\cdots,N-1)$ be the number of entries on the $a$-th skew-diagonal of  an $L\times K$ matrix, 
\begin{align*}
w_a = \#\lcb (i,j)~|~ i+j=a,~ 0\leq i\leq L-1,~0\leq j\leq K-1\rcb,\numberthis\label{eq:weight}
\end{align*}
the $a$-th entry of $\H^\dagger\BZ$, denoted $[\H^\dagger\BZ]_a$  is given by 
\begin{align*}
[\H^\dagger\BZ]_a=\frac{1}{w_a} \sum_{i+j=a}[\BZ]_{ij}.
\end{align*}
For example, consider the following $3\times 3$ matrix,
\begin{align*}
\BZ=\begin{bmatrix}
1 & 2 & 3\\
4 & 5 & 6\\
7 & 8 & 9
\end{bmatrix},
\end{align*}
we have 
\begin{align*}
w_0=1,~w_1=2,~w_2=3,~w_3=2,~w_4=1,
\end{align*}
and
\begin{align*}
[\H^\dagger\BZ]_0=1,~[\H^\dagger\BZ]_1=\frac{2+4}{2},~[\H^\dagger\BZ]_2=\frac{3+5+7}{3},~[\H^\dagger\BZ]_3=\frac{6+8}{2},~[\H^\dagger\BZ]_4=9.
\end{align*}
Moreover, it is not hard to show that $\H^\dagger\H=\I$, where $\I$ is the identity operator from $\C^N$ to $\C^N$.

We are now in position to introduce the Cadzow's algorithm which is based on the fact that signals of interest arising from a wide range of applications have low rank   Hankelization. That is, for the denoising problem \eqref{eq:denoising} we can  assume $\rank(\H\bx)=r$, where $r\ll \min(L,K)$. Additive noise often increases the rank of $\H\bx$ and hence it is very typical that $\rank(\H\by)=\min(L,K)$. Consequently, an intuitive way of estimating $\bx$ from $\by$ is to do rank reduction. Starting from $\bz_0=\by$, the Cadzow's algorithm iteratively updates the estimate via the following rule:
\begin{align*}
\bz_{k+1} = \H^\dagger\T_r\H\bz_k,~k=0,\cdots.\numberthis\label{eq:cadzow}
\end{align*}
Here $\T_r$ computes the truncated SVD of an $L\times K$ matrix, that is,
\begin{align*}
\mathcal{T}_r(\bm{Z})=\sum_{j=1}^r\sigma_j\bm{u}_j\bm{v}_j^*,
\quad\mbox{where }\bm{Z}=\sum_{j=1}^{\min(L,K)}\sigma_j\bm{u}_j\bm{v}_j^*\mbox{ is an SVD with }\sigma_1\geq\sigma_2\geq\ldots\geq\sigma_{\min(L,K)}.
\end{align*}
It is worth noting that the Cadzow's algorithm is also known as multichannel
singular spectrum analysis (MSSA) which was proposed  by Broomhead and King in 1986 for the analysis of time series from dynamical system \cite{broomhead1986extracting}. In particular, the first iteration of the Cadzow's algorithm is usually referred to as singular spectrum analysis (SSA) which is also very important in many applications; see for example \cite{hassani2010review,hassani2009singular,golyandina2001analysis,golyandina2005ssa}.

We can also interpret the Cadzow's algorithm as the method of alternating projections in the matrix domain. To this end, let $\M_r$ and $\M_{H}$ be the set of rank-$r$ and Hankel matrices, respectively. For any matrix $\BZ\in\C^{L\times K}$, it follows from the Eckart-Young theorem that the projection of $\BZ$ onto $\M_r$ is given by $\T_r(\BZ)$. Moreover, it can be easily verified that the projection of $\BZ$ onto $\M_H$ is given by $\H\H^\dagger\BZ$. Therefore, if we let $\BZ_k=\H\bz_k$, after multiplying $\H$ on both sides of \eqref{eq:cadzow}, it can be seen that the Cadzow's algorithm is equivalent to 
\begin{align*}
\BZ_{k+1} = \P_{\M_H}\P_{\M_r}\BZ_k,~k=0,\cdots,\numberthis\label{eq:cadzow_matrix}
\end{align*}
where $\P_{\M_r}$ and $ \P_{\M_H}$  denote the projection onto $\M_r$ and $\M_H$, respectively.

As suggested in \cite{oropeza2011simultaneous}, the Cadzow's algorithm can be adapted for signal recovery problems when there are missing entries. Recall that $\bx$ is our target signal. Suppose in some situations we are not able to observe all  the entries of $\bx$, but can only observe those entries with indices  in $\Omega$, where $\Omega$ is a subset of $\{0,\cdots,N-1\}$. Moreover, let $\P_\Omega(\bx)$ denotes the samples of $\bx$ on $\Omega$, namely, 
\begin{align*}
[\P_\Omega(\bx)]_{j} =\begin{cases}
x_j & \mbox{if }j\in\Omega\\
0 & \mbox{otherwise}.
\end{cases}
\end{align*}
A natural question is whether it is possible to reconstruct $\bx$ from $\P_\Omega(\bx)$. This is an ill-posed problem without any restrictions on $\bx$. However, when $\H\bx$ is low rank, it has been shown that the missing entries can be completed through convex methods \cite{chen2014robust} as well as non-convex methods \cite{cai2019fast}. In \cite{oropeza2011simultaneous}, a variant of the Cadzow's algorithm was proposed for signal recovery. The algorithm iteratively fills the missing entries in the following way: $\bz_0=\P_\Omega(\bx)$,
\begin{align*}
\bz_{k+1} = \P_\Omega(\bx)+(\I-\P_\Omega)\H^\dagger\T_r\H\bz_k,~k=0,\cdots.\numberthis\label{eq:cadzow_missing}
\end{align*}
Roughly speaking, the algorithm refines the estimation of the unknown entries using the entries of $\H^\dagger\T_r\H\bz_k$  on $\Omega^c$. In addition, when noise exists on the observed entries, i.e., under the measurement model $\P_\Omega(\by)=\P_\Omega(\bx)+\P_\Omega(\be)$, one can use an appropriate linear combination of $\P_\Omega(\by)$ and $\P_\Omega(\H^\dagger\T_r\H\bz_k)$ to estimate $\P_\Omega(\bx)$ rather than simply filling back the noisy entries. Interested readers are referred to \cite{oropeza2011simultaneous} for details.
\subsection{Examples where low rank Hankelization appears}\label{subsec:intro1}
As shown previously,  the low rank structure of Hankel matrices plays a key role in the development of the Cadzow's algorithm and its variant. In this subsection we present several  real examples  
where low rank Hankelization appears. 

\paragraph{Time series satisfying an LRF} Let $\bx=[x_0,\cdots,x_{N-1}]^\top$ be a time series of length $N$. We say $\bx$ satisfies a linear recurrent formula (LRF) if 
\begin{align*}
x_j = a_1x_{j-1}+\cdots +a_rx_{j-r},\quad j=r,\cdots, N-1,\numberthis\label{eq:LRF}
\end{align*}
for some coefficients $a_1,\cdots, a_r$. Let $L<N$ be an integer which is usually referred to as the window length, and let $K=N-L+1$. Define the $L$-lagged vectors $\bx_j = [x_j,\cdots,x_{j+L-1}]^\top,~j = 0,\cdots,K-1$. If we construct an $L\times K$ matrix $\BX$ as follows 
\begin{align*}
\BX = [\bx_0,\cdots,\bx_{K-1}], 
\end{align*}
it is easy to see that $\BX$ is indeed a Hankel matrix and in fact we have $\H\bx=\BX$.
Moreover, since $\bx$ is a time series which obeys \eqref{eq:LRF}, $\H\bx$ is a matrix of rank at most $r$. Time series satisfying an LRF is a very important model for many real applications, see \cite{golyandina2001analysis} for more details. 

\paragraph{Spectrally sparse signals} Consider the 1D spectrally sparse signal $x(t)$ consisting of $r$ 
complex sinusoids  without damping factors
\begin{equation}\label{eq:spectral_sparse}
x(t)=\sum_{j=1}^{r}d_je^{ i2\pi f_jt},
\end{equation}
where $f_j\in[0,1)$ is the normalized frequency, and $d_k\in\mathbb{C}$ is the corresponding complex amplitude. Let $\bx=[x(0),\cdots,x(N-1)]^\top$ be the discrete samples of $x(t)$ at $t\in\{0,\cdots,N-1\}$.
Letting $y_j=e^{i2\pi f_j}$, because $\bx$ is obtained from a spectrally sparse signal, it can be shown that \cite{hua1992rank,yang1996rank} $\H\bx$ admits the Vandermonde  decomposition of the form
$$
\mathcal{H}\bm{x}=\bm{E}_L\bm{D}\bm{E}_R^T,
$$
where
$$
\bm{E}_L=
\left[\begin{array}{cccc}
1         & 1         & \cdots & 1 \\
y_1       & y_2       & \cdots & y_r \\
\vdots    & \vdots    & \vdots & \vdots \\
y_1^{L-1} & y_2^{L-1} & \cdots & y_r^{L-1} \\
\end{array}\right],~
\bm{E}_R=
\left[\begin{array}{cccc}
1         & 1         & \cdots & 1 \\
y_1       & y_2       & \cdots & y_r \\
\vdots    & \vdots    & \vdots & \vdots \\
y_1^{K-1} & y_2^{K-1} & \cdots & y_r^{K-1} \\
\end{array}\right]
$$
and
$\bm{D}$ is a diagonal matrix whose diagonal entries are  $d_1,\dots,d_r$. The Vandermonde decomposition of $\H\bx$ implies  that $\rank(\H\bx)=r$ provided $d_j\neq 0$ and $f_j\neq f_{j'}$ when $j\neq j'$.

\paragraph{Periodic stream of Diracs} In the previous two examples, low rank Hankelization occurs in the time domain. There are also examples where low rank Hankelization occurs in the Frequency domain. Consider the following $1$-periodic of Diracs studied in \cite{blu2008sparse}:
\begin{align*}
x(t) = \sum_{j=1}^r \sum_{j'\in\mathbb{Z}}x_j\delta(t-t_j-j'),\numberthis\label{eq:dirac_signal}
\end{align*}
where $x_j\in\R$ and $t_j\in[0,1)$. Since $x(t)$ is $1$-periodic,  it can be represented using the Fourier series; that is
\begin{align*}
x(t) = \sum_{k\in\mathbb{Z}}\hat{x}_k e^{i 2\pi kt},\quad \mbox{where }\hat{x}_k=\sum_{j=1}^rx_j e^{-i 2\pi kt_j}.
\end{align*}
Let $\hat{\bx}$ be a vector consisting of $N$ consecutive values of $\hat{x}_k$ and define 
\begin{align*}
\hat{x}(\omega) = \sum_{j=1}^rx_j e^{-i 2\pi t_j\omega}.\numberthis\label{eq:time_sparse}
\end{align*}
We can view $\hat{\bx}$ as a discrete sample vector from $\hat{x}(\omega)$. Noticing the similarity between \eqref{eq:spectral_sparse} and \eqref{eq:time_sparse}, it is not hard to see that  $\rank(\H\hat{\bx})=r$. 

\subsection{Main contributions and outline}
The main contributions of this paper are two-fold. Firstly, a Fast Cadzow's algorithm has been developed through the introduction of an additional subspace projection; see Section~\ref{sec:algorithm}. Numerical experiments in Section~\ref{sec:numerics} establish the computational advantages of the Fast Cadzow's algorithm over the Cadzow's algorithm. Secondly, to address the non-decreasing MSE issue of the Cadzow's and Fast Cadzow's algorithms for signal denoising, a Gradient method and its accelerated variant are introduced in Section~\ref{sec:gradient}. The algorithms are then tested for different denoising problems to justify the effectiveness of the gradient algorithms. Section~\ref{sec:conclusion} concludes this paper with a short discussion.

\section{Fast Cadzow's algorithm}\label{sec:algorithm}
In general, the complexity of computing the SVD of an $L\times K$ $(L+K-1=N)$ matrix is $O(N^3)$ when $L$ and $K$ are proportional to $N$, causing the computation of the SVD (i.e., the $\T_r$ operation) to be the dominant cost in the Cadzow's iteration \eqref{eq:cadzow}. This has limited the computational efficiency of the Cadzow's algorithm, especially for high dimensional problems. In this section, we  present an accelerated variant of the Cadzow's algorithm, dubbed Fast Cadzow's algorithm, for  the   denoising  problem as well as the  recovery problem from missing entires; see  Algorithm~\ref{alg:fastCadzow} for a formal description. We only focus on the discussion of the algorithm for the denoising problem since the same idea is used for the recovery problem.

\begin{algorithm}[ht!]
\caption{Fast Cadzow's algorithm}
\label{alg:fastCadzow}
\begin{algorithmic}[1]
\Require $\bz_0=\by$
\For{$k=0,1,\cdots$}
\State Compute a matrix subspace $T_k$
\If{ \{for signal denoising\}}
\State $\bz_{k+1}=\H^\dagger\T_r\P_{T_k}\H\bz_k$
\Else { \{for recovery from missing entries\}}
\State $\bz_{k+1} = \P_\Omega(\bx)+(\I-\P_\Omega)\H^\dagger\T_r\P_{T_k}\H\bz_k$
\EndIf
\EndFor
\Ensure $\bz_k$
\end{algorithmic}
\end{algorithm}

\begin{figure}
\centering
\includegraphics[height=5cm,width=7cm]{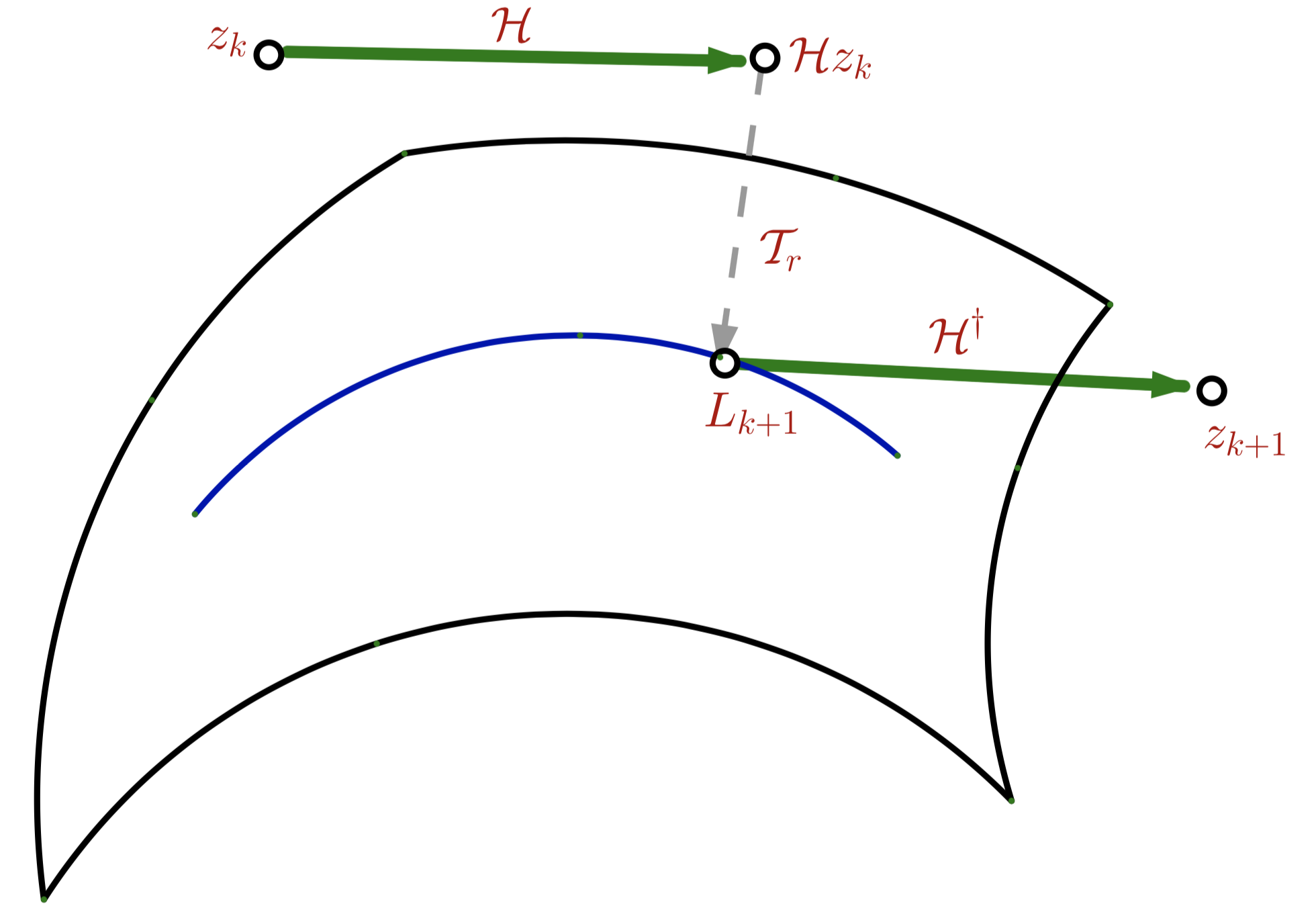}\hspace{0.5cm}
\includegraphics[height=5cm,width=7cm]{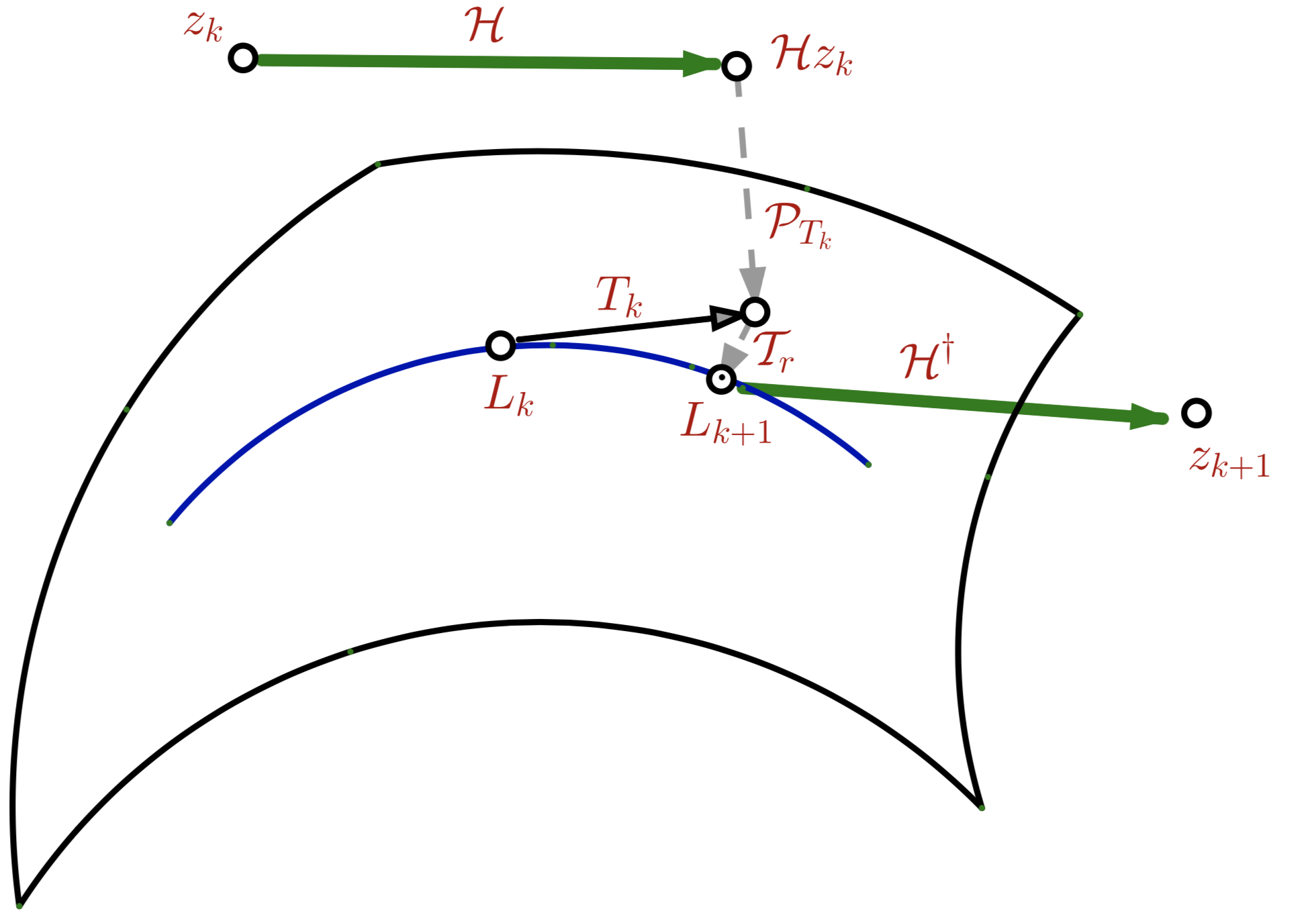}
\caption{Pictorial illustrations of Cadzow's algorithm (left) and Fast Cadzow's algorithm (right).}\label{fig:algs_pic}
\end{figure}

Recall that the denoising problem is about estimating a signal $\bx$ from the noisy measurements $\by=\bx+\be$. We assume the Hankel matrix $\H\bx$ associated with $\bx$ is rank $r$. The Fast Cadzow's iteration is overall similar to the Cadzow's iteration,  except that there is an additional projection $\P_{T_k}$ in the Fast Cadzow's algorithm, where $\P_{T_k}$ denotes the operator of projecting a matrix onto a subspace $T_k$; see \eqref{eq:projection} on how to compute $\P_{T_k}$ for a particular choice of $T_k$. In other words, after the Hankel matrix $\H\bz_k$ is constructed, we first project it onto $T_k$, followed by the truncation to the best rank $r$ approximation. The motivation behind Algorithm~\ref{alg:fastCadzow} is that after $\P_{T_k}$ is introduced the resulting matrix can be more structured. Thus it can be expected that the SVD of the matrix can be computed in a fast way. Notice that if we choose $T_k$ to be $\C^{L\times K}$ in each iteration, Algorithm~\ref{alg:fastCadzow} is indeed the Cadzow's algorithm; see Figure~\ref{fig:algs_pic} for pictorial illustrations of the Cadzow's algorithm  and the Fast Cadzow's algorithm. Next, we will see how to choose the subspace $T_k$ in a clever way such that the truncation operator $\T_r$ can be performed efficiently.
 
\subsection{Choice for $T_k$}
Our choice of $T_k$ is inspired by  the Riemannian optimization methods for low rank matrix recovery \cite{WCCL:SIMAX:16,vandereycken2013low}, which is closely related to the manifold structure of low rank matrices. To motivate this choice,  assume  $\bz_k$ has already been computed via $\bz_{k}=\H^\dagger\T_r\P_{T_{k-1}}\H\bz_{k-1}$ in the $(k-1)$-th iteration. Let $\BL_k = \T_r\P_{T_{k-1}}\H\bz_{k-1}$. Noticing that $\BL_k$ is a matrix of rank $r$, it has the reduced SVD  $\BL_k=\BU_k\BS\BV_k^*$, where $\BU_k$ and $\BV_k$ are $L\times r$ and $K\times r$ orthogonal matrices, respectively. In the $k$-th iteration, $T_k$ is selected to be the direct sum of the column and row subspaces of $\BL_k$, i.e,. 
\begin{align*}
T_k = \{\BU_k\BB^*+\BC\BV_k^*~|~\BB\in\C^{K\times r},~\BC\in\C^{L\times r}\}.\numberthis\label{eq:tangentSpace}
\end{align*}
From the perspective of differential geometry, the set of fixed rank $r$ matrices forms a smooth manifold and $T_k$ is  the tangent space of the manifold at $\BL_k$. 


Given a matrix $\BZ\in\C^{L\times K}$, the projection of $\BZ$ onto $T_k$ is given by 
\begin{align*}
\P_{T_k}(\BZ) = \BU_k\BU_k^*\BZ+\BZ\BV_k\BV_k^*-\BU_k\BU_k^*\BZ\BV_k\BV_k^*.
\numberthis\label{eq:projection}
\end{align*}
 In the first iteration (i.e., $k=0$) of Algorithm~\ref{alg:fastCadzow}, we simply choose $T_0\in\C^{L\times K}$.  That is, the first iteration of the Fast Cadzow's algorithm  coincides with that of the Cadzow's algorithm, which is  an SSA step.
Noting that $\bz_k=\H^\dagger\BL_k$ and $\H\H^\dagger\BL_k\approx \BL_k$ when $\BL_k$ is close to a Hankel matrix,  we have $\H\bz_k\approx \BL_k$ in this situation. It follows that $\P_{T_k}\H\bz_k\approx \H\bz_k$ since $\BL_k\in T_k$. Therefore, the projection of $\H\bz_k$ onto $T_k$ can capture most of its energy and we can expect that the Fast Cadzow's algorithm should exhibit similar behavior to the Cadzow's algorithm. The numerical results in Section~\ref{sec:numerics} will confirm this intuition. 
\subsection{Computational complexity}
The true novelty of introducing the additional projection $\P_{T_k}$ is that after this projection the SVD of the matrix can be computed more efficiently. In a nutshell, the SVD of an $L\times K$ matrix can be reduced to the SVD of a $2r\times 2r$ matrix.
In this subsection we will investigate the computational complexity of the Fast Cazdow's algorithm, together with the details for computing the SVD. 

For ease of exposition we split the single Fast Cadzow's iteration into three steps:
\begin{align*}
\begin{cases}
\BW_k = \P_{T_k}\H\bz_k\\
\BL_{k+1} = \T_r\BW_k\\
\bz_{k+1} =\H^\dagger \BL_{k+1}.
\end{cases}\numberthis\label{eq:algSplit}
\end{align*}
Letting $\BH_k=\H\bz_k$ which will not be formed explicitly, 
by \eqref{eq:projection}, we have 
\begin{align*}
\BW_k&=\BU_k\BU_k^*\BH_k+\BH_k\BV_k\BV_k^*-\BU_k\BU_k^*\BH_k\BV_k\BV_k^*\\
&=\BU_k\BU_k^*\BH_k\BV_k\BV_k^*+\BU_k\BU_k^*\BH_k(\BI-\BV_k\BV_k^*)+(\BI-\BU_k\BU_k^*)\BH_k\BV_k\BV_k^*\\
&=\BU_k\BG\BV_k^* + \BU_k\BB^*+\BC\BV_k^*,\numberthis\label{eq:W}
\end{align*}
where $\BG=\BU_k^*\BH_k\BV_k$, $\BB=(\BI-\BV_k\BV_k^*)\BH_k^*\BU_k$ and $\BC=(\BI-\BU_k\BU_k^*)\BH_k\BV_k$. In a practical implementation $\BW_k$ will be stored in the form of \eqref{eq:W}. Since $\BH_k$ is a Hankel matrix, both $\BH_k\BV_k$ and $\BH_k^*\BU_k$ can be computed using $r$ fast matrix-vector multiplications without forming $\BH_k$. This costs $O(Nr\log N)$ flops assuming $L\sim K\sim N$. Once $\BH_k\BV_k$ and $\BH_k^*\BU_k$ are known, $\BG$, $\BB$ and $\BC$ can be computed using a few matrix-matrix products with $O(Nr^2)$ flops. Thus, it requires $O(Nr^2+Nr\log N)$ flops to compute $\BG$, $\BB$ and $\BC$ in \eqref{eq:W}.

Next we show how to reduce the SVD of $\BW_k$ to the SVD of a $2r\times 2r$ matrix starting from the decomposition \eqref{eq:W}.  Let $\BB= \BQ_1\BR_1$ and $\BC=\BQ_2\BR_2$ be the QR decompositions of $\BB$ and $\BC$, respectively. Then, it is not hard to see that 
\begin{align*}
\BQ_1\perp\BV_k,\quad\BQ_2\perp\BU_k.
\end{align*}
Moreover, we have 
\begin{align*}
\BW_k& = \BU_k\BG\BV_k^* + \BU_k\BR_1^*\BQ_1^*+\BQ_2\BR_2\BV_k^*\\
&=\begin{bmatrix}\BU_k & \BQ_2\end{bmatrix}
\begin{bmatrix}
\BG& \BR_1^*\\
\BR_2 & \bm{0}
\end{bmatrix}
\begin{bmatrix}
\BV_k&\BQ_1
\end{bmatrix}^*.\numberthis\label{eq:decW}
\end{align*}
Let $\BM$ be the middle $2r\times 2r$ matrix and suppose its SVD is given by 
\begin{align*}
\begin{bmatrix}
\BG& \BR_1^*\\
\BR_2 & \bm{0}
\end{bmatrix} = {\BU}_G{\BS}_G{\BV}_G^*.
\end{align*}
Since both $\begin{bmatrix}
\BU_k & \BQ_2
\end{bmatrix}$ and $\begin{bmatrix}
\BV_k&\BQ_1
\end{bmatrix}$ are orthogonal matrices, the SVD of $\BW_k$ is given by 
\begin{align*}
\BW_k = \lb\begin{bmatrix}\BU_k & \BQ_2\end{bmatrix}\BU_G\rb\BS_G\lb\begin{bmatrix}
\BV_k&\BQ_1
\end{bmatrix}\BV_G\rb^*.
\end{align*}
From the above discussion, we can see that computing $\BL_{k+1}$ from $\BW_k$ (or equivalently, computing the SVD of $\BW_k$) requires $O(Nr^2+r^3)$ flops which account for the QR decompositions of $\BB$ and $\BC$, and the SVD  of  $\BM$. Moreover, when $L=K$ (i.e., matrices are square), nearly half computational costs and storage can be saved by using the Takagi factorization. Interested readers can find the details in \cite{cai2019fast}.

It remains to investigate the cost for computing $\bz_{k+1}=\H^\dagger\BL_{k+1}$. Let
$\BL_{k+1}=\BU_{k+1}\BS_{k+1}\BV_{k+1}^*$ be the SVD of $\BL_{k+1}$ which can be  obtained by truncating the SVD of $\BW_k$. We have 
\begin{align*}
\H^\dagger\BL_{k+1} =\sum_{j=1}^r\BS_{k+1}(j,j)\H^\dagger\lsb\BU_{k+1}(:,j)\lb\BV_{k+1}(:,j)\rb^*\rsb.
\end{align*}
Noting that 
\begin{align*}
\lsb\H^\dagger\lsb\BU_{k+1}(:,j)\lb\BV_{k+1}(:,j)\rb^*\rsb\rsb_a =\frac{1}{w_a}\sum_{p+q=a} \BU_{k+1}(p,j)\overline{\BV}_{k+1}(q,j),\quad a=0,\cdots, N-1,
\end{align*}
where $w_a$ is defined in \eqref{eq:weight}, $\H^\dagger\BL_{k+1}$ can be computed by fast convolution using $O(N\log N)$ flops. Thus, computing $\bz_{k+1}$ from $\BL_{k+1}$ requires $O(Nr\log N)$ flops.

Putting it all together, the dominant per iteration computational cost of the Fast Cadzow's algorithm is $O(Nr^2+Nr\log N+r^3)$. In addition, $O(Nr)$ space is required
to store the matrices appearing in the SVD and the QR decompositions. 
\subsection{High dimensional problems}
We have presented the Cadzow's and Fast Cadzow's algorithms for the 1D problem.
However, the algorithms can also work for high dimensional problems based on the multi-fold Hankel structures. For ease of exposition we discuss the two-dimensional case briefly  but emphasize that the situation in general higher dimensional cases is similar.

Let $\BX\in\R^{N_1\times N_2}$ be a 2D signal, and let $(L_1,K_1)$ and $(L_2,K_2)$ be two pairs of positive numbers which satisfy $L_1+K_1-1=N_1$ and $L_2+K_2-1=N_2$. The Hankel matrix corresponding to $\BX$, denoted $\H\BX$, is defined as follows: 
\begin{align*}
\H\BX=\begin{bmatrix}
\H\BX(0,:) & \H\BX(1,:) & \H\BX(2,:) & \cdots & \H\BX(K_1-1,:)\\
\H\BX(1,:) & \H\BX(2,:) & \H\BX(3,:) & \cdots & \H\BX(K_1,:)\\
\H\BX(2,:) & \H\BX(3,:) & \H\BX(4,:) & \cdots & \H\BX(K_1+1,:)\\
\vdots & \vdots & \vdots & \vdots & \vdots\\
\H\BX(L_1-1,:) & \H\BX(L_1,:)& \H\BX(L_1+1,:) & \cdots &\H\BX(N_1-1,:)
\end{bmatrix},
\end{align*}
where $\H X(i,:)\in\R^{L_2\times K_2}$ is the Hankel matrix for the $i$-th row of $\BX$.
In other words, we first form the Hankel matrix for each row of $\BX$ and then form the block Hankel matrix $\H\BX$ using  the row Hankel matrices. 

Assuming $\H\BX$ is a low rank matrix, the Cadzow's and Fast Cadzow's algorithms can be similarly developed for the 2D denoising and reconstruction problems. The details will be omitted. Moreover, there are real signals whose Hankel matrices are low rank, for example if $\BX$ is a 2D spectrally sparse signal or a frequency slice of seismic data after the Fourier transform \cite{hua1992rank,yang1996rank,oropeza2010singular}.

\section{Numerical Experiments}\label{sec:numerics}
In this section we evaluate the empirical performance of the Fast Cadzow's algorithm against the Cadzow's algorithm. In our implementations of the algorithms $L$ and $K$ (or $L_i$ and $K_i$ for higher dimensions) are chosen in such a way that the resulting Hankel matrix is approximately square.
 Numerical results demonstrate that the two algorithms exhibit similar denoising and reconstruction performance while the Fast Cadzow's algorithm can be substantially faster. 
 We compare the two algorithms on three different settings:
 \begin{itemize}
 \item in Subsection~\ref{subsec:numerics1}, the algorithms are tested on randomly generated spectrally sparse signals;
 \item in Subsection~\ref{subsec:numerics2}, the algorithms are tested on periodic stream of Diracs in which the low rank structure is hidden in the Fourier domain;
 \item in Subsection~\ref{subsec:numerics3}, the algorithms are tested on seismic data denoising and reconstruction.
 \end{itemize}
 The experiments are executed from Matlab 2018a on a desktop with two Inter i5 CPUs (1.60GHz and 1.80GHz respectively) and 8GB memory. 

To utilize the fast Hankel matrix-vector multiplication, we use the Krylov subspace-based SVD solver  from the PROPACK package \cite{PROPACK} to compute the partial SVD of $\H\bz_k$ in the Cadzow's algorithm. It would be difficult to summarize the computational cost of such a method for computing the partial SVD and optimistically we can say that the cost is $O(Nr^2+Nr\log N)$ flops.  Though it has the same order as the SVD method for the Fast Cadzow's algorithm, it is worth noting that the constant hidden in  $O(Nr^2+Nr\log N)$ is not clear since the performance of the Krylov's method {\em depends heavily on properties of the input matrix and on the amount of
effort spent to stabilize the algorithm} \cite{Halko2011randomizedSVD}. In contrast, the constants in  the costs $O(Nr^2)$ and $O(Nr\log N)$ for computing the SVD in the Fast Cadzow's algorithm are exactly known since for example $O(Nr\log N)$ only comes from the computations of $\BH_k\BV_k$, $\BH_k^*\BU_k$ and $\H^\dagger\BL_k$. Moreover, our numerical experiments  show that the Fast Cadzow's algorithm is still two to four times faster  even after the fast partial  SVD package is used for the Cadzow's algorithm.

\paragraph{Remark on the application of randomized SVD in Cadzow's algorithm} We may also use randomized SVD  \cite{RST09,Halko2011randomizedSVD} to compute the low rank approximation of $\H\bz_k$ in the Cadzow's algorithm. In fact, both randomized SVD and the accelerated method used in the Fast Cadzow's algorithm are based on the subspace projection. The two methods differ mainly in the manner they choose the subspace. The randomized SVD computes a subspace in each iteration randomly while the method used in the Fast Cadzow's algorithm exploits the geometric structure of the low rank matrix manifold and utilizes the tangent space in each iteration. Letting $\BH_k=\H\bz_k$, consider the following randomized SVD procedure  in \cite{Halko2011randomizedSVD} for computing a good rank-$r$ approximation of $\H\bz_k$:
\begin{enumerate}
\item Generate an $L\times 2r$ random Gaussian matrix $\BOMG$  and compute $\BZ=(\BH_k\BH_k^*)^q\BH_k\BOMG$.
\item Compute the reduced QR factorization of $\BZ$: $[\BQ,\BR]=\mbox{qr}(\BZ)$.
\item Compute the reduced SVD of $\BB=\BQ^*\BA$: $[\BU_B,\BS_B,\BV_B]=\mbox{svd}(\BB)$.
\item Construct the rank-$r$ approximation of  $\BH_k$ as follows: $\BH_k\approx[\BU_k,\BS_k,\BV_k]$ with $\BU_k=\BQ\BU_B$, $\BS_k=\BS_B$, and $\BV_k=\BV_B$.
\end{enumerate}

The matrix power $(\BH_k\BH_k^*)^q$ in the first step is  introduced to accelerate the decay of the singular spectrum of the matrix, and typical choices for $q$ are $q=1$ or $q=2$ in practice. By counting the number of flops, it is not hard to see that the cost of computing the left and right matrices in \eqref{eq:decW} is overall similar to that of computing $\BH_k\BOMG$ and $\BB$ in the above randomized SVD procedure. However, due to the extra cost incurred by the matrix power $(\BH_k\BH_k^*)^q$ and the fact that $\BB$ is an $2r\times K$ matrix while the middle matrix in \eqref{eq:decW} is only a $2r\times 2r$ matrix, the total computational cost of the randomized SVD is higher than that of the accelerated method in the Fast Cadzow's algorithm. Moreover, since $\BB$ becomes an unstructured matrix after projection, we cannot apply the fast matrix-vector product when computing its partial SVD. Indeed, our numerical tests show that the Cadzow's algorithm with the randomized SVD is not as efficient as the Cadzow's using the PROPACK package for the computation of the partial SVD. Here we will omit the computational results of the Cadzow's algorithm with the randomized SVD.
\subsection{Spectrally sparse signal denoising and reconstruction}\label{subsec:numerics1}

The definition of spectrally sparse signals can be found in Section~\ref{subsec:intro1}. We first test the performance of the algorithms on signal denoising. The test signals are generated in the following way: the frequency of each harmonic is uniformly sampled from $[0,1)$, and the argument of the weight is uniformly sampled from $[0,2\pi)$ while the amplitude is chosen to be $1+10^{0.5c}$ with $c$ being uniformly sampled from $[0,1]$. 

\begin{table}[ht!]
\centering
\caption{Average computational time (seconds), average number of iterations and average mean square error of Cadzow and Fast Cadzow over 10 random problem instances per $(N,r)$ for $N\in\lbrace 4096, 128\times 128, 64\times64\times64\rbrace$ and  $r \in \lbrace 5,10,20\rbrace$.
}
\label{tab:spectral_denoise}
\vspace{0.2cm}
\begin{tabular}{|c|ccc|ccc|ccc|}
\hline
$r $         & \multicolumn{3}{c|}{5}    & \multicolumn{3}{c|}{10}  &\multicolumn{3}{c|}{20}  \\
\hline           & MSE    & Iter   & Time    & MSE    & Iter   & Time  & MSE    & Iter   & Time  \\
\hline           & \multicolumn{9}{c|}{$N = 4096$}                            \\
\hline Cadzow     & 3.16e-02 & 10.3 & 0.48  & 4.17e-02 & 11.1 & 0.80  & 6.15e-02 & 11.4 & 1.75\\
\hline Fast Cadzow & 3.16e-02 & 10.3 & 0.19 & 4.17e-02 & 11.1 & 0.38 & 6.15e-02 & 11.4 & 0.81 \\
\hline         
& \multicolumn{9}{c|}{$N = 128\times128$}                                                         \\
\hline Cadzow               & 2.02e-02 & 14.6 & 0.61  & 2.95e-02 & 15.1 & 1.09  & 4.19e-02 & 15.8 & 2.98  \\
\hline Fast Cadzow           & 2.02e-02 & 14.6 & 0.28  & 2.95e-02 & 15.1 & 0.47  & 4.19e-02 & 15.8 & 1.04      \\
\hline         & \multicolumn{9}{c|}{$N = 64\times64\times64$}                                                          \\
\hline Cadzow               & 0.70e-02 & 13.1 & 10.45  & 0.96e-02 & 13.8 & 18.18  & 1.35e-02 & 14.0 & 40.98  \\
\hline Fast Cadzow           & 0.70e-02 & 13.1 & 4.61  & 0.96e-02 & 13.8 & 9.66 & 1.35e-02 & 14.0 & 20.25  \\
\hline   
\end{tabular}
\end{table}

We consider the additive noise model in which a true signal $\bx$ is corrupted by a noisy vector $\be$ in the form of
\begin{align*}
\be = \epsilon\cdot\|\bx\|\cdot\frac{\bw}{\|\bw\|},
\end{align*}
where $\epsilon$ is referred to the noise level and $\bw$ follows the standard multivariate normal distribution. Letting $\bz_k$ be the output of the algorithms,  the mean squared error (MSE) defined by 
\begin{align*}
\mbox{MSE} = \frac{\|\bz_k-\bx\|}{\|\bx\|}
\end{align*}
will be used to evaluate their performance.
Tests are conducted for 1D, 2D and 3D spectrally sparse signals with $N\in\{ 4096, 128\times 128, 64\times64\times64\}$, $r\in \{5,10,20\}$,  and $\epsilon = 0.5$. The Cadzow's and Fast Cadzow's algorithms are terminated whenever $\|\bz_{k+1}-\bz_k\|/\|\bz_k\|\leq 10^{-6}$. For each pair of $(N,r)$ $10$ randomly generated problem instances are tested. Then we compute the average  computation time, the average number of iterations and the average MSE for each algorithm; see Table~\ref{tab:spectral_denoise}. From the table we can see that the Cadzow's and Fast Cadzow's algorithms achieve similar MSE and it also takes them roughly the same number of iterations  to converge. However, the Fast Cadzow's algorithm is at least two times faster than the Cadzow's algorithm since the additional subspace projection can reduce the computational cost of the SVD.

\begin{table}[ht!]
\caption{Average computational time (seconds), average number of iterations and average mean square error of Cadzow and Fast Cadzow over 10 random problem instances per $(n,r)$ for $N\in\lbrace 4096, 128\times 128, 64\times64\times64\rbrace$ and  $r \in \lbrace 5,10,20\rbrace$.
}
\centering
\label{tab:spectral_missing}
\vspace{0.2cm}
\begin{tabular}{|c|ccc|ccc|ccc|}
\hline
$r$          & \multicolumn{3}{c|}{5}    & \multicolumn{3}{c|}{10}  &\multicolumn{3}{c|}{20}  \\
\hline           & MSE    & Iter   & Time    & MSE   & Iter  & Time  & MSE    & Iter  & Time  \\
\hline           & \multicolumn{9}{c|}{$N = 4096$}                            \\
\hline Cadzow     & 8.47e-11 & 35.9 & 1.85  & 9.09e-11 & 38.0 & 3.08  & 1.10e-10 & 41.1 & 8.43\\
\hline Fast Cadzow & 8.47e-11 & 35.9 & 0.61  & 9.09e-11 & 38.0 & 1.35 & 1.10e-10 & 41.1 & 3.45\\
\hline           & \multicolumn{9}{c|}{$N = 128\times128$}                            \\
\hline Cadzow     & 7.45e-11 & 35.3 & 1.50  & 8.63e-11 & 36.4 & 2.73 & 9.63e-11 &38.2 & 8.00 \\
\hline Fast Cadzow & 7.45e-11 & 35.3 & 0.63  & 8.63e-11 & 36.4 & 1.09 & 9.63e-11 &38.2 & 2.57\\
\hline                      & \multicolumn{9}{c|}{$N = 64\times64\times64$}                                                         \\
\hline Cadzow               & 7.12e-11 &34.0 & 30.08  & 7.55e-11 & 34.1 & 47.98  & 9.02e-11 & 34.2 & 108.86  \\
\hline Fast Cadzow           & 7.12e-11 & 34.0 & 11.61  & 7.55e-11 & 34.1 & 22.76  & 9.02e-11 & 34.2 & 47.28  \\ \hline 
\end{tabular}
\end{table}

Next we compare the two algorithms on reconstruction problems when there are missing entries. The overall setup is similar to the denoising case except that instead of 
having $\bx$ being contaminated by additive noise we only observe $50\%$ of its entries uniformly at random. Note that the variant of the Cadzow's algorithm for handling missing entries is presented in \eqref{eq:cadzow_missing}. In the tests the algorithms are terminated when $\|\bz_{k+1}-\bz_k\|/\|\bz_k\|\leq 10^{-10}$. The average computational results over $10$ random simulations are presented in Table~\ref{tab:spectral_missing}. Clearly, both of the algorithms are able to exactly recover the missing entries and the Fast Cadzow's algorithm is substantially faster.

The  algorithms can also handle missing entries and additive noise simultaneously. The basic idea is to make an appropriate linear combination between observed measurements and the update. More precisely, we consider the following variants  of the algorithms: 
\begin{align*}
     \mbox{(Cadzow)}~ & \bz_{k+1} = \alpha\cdot\mathcal{P}_{\Omega}(\by)+(\mathcal{I}-\alpha\cdot\mathcal{P}_{\Omega})\H^\dagger\T_r\H\bz_k,\numberthis\label{eq:cadzow_denoising_missing}\\
     \mbox{(Fast Cadzow)}~ & \bz_{k+1} = \alpha\cdot\P_\Omega(\by)+(\I-\alpha\cdot\P_\Omega)\H^\dagger\T_r\P_{T_k}\H\bz_k,\numberthis\label{eq:fastcadzow_denoising_missing}
\end{align*}
where $\alpha\geq 0$ is a tuning parameter. In our tests we choose $\alpha =0.8$ as suggested by Gao et al. \cite{gao2013fast}, and the algorithms are terminated when $\|\bz_{k+1}-\bz_k\|/\|\bz_k\|\leq 10^{-6}$. The numerical results for this setting are summarized in Table~\ref{tab:spectral_missing_noise}. Again, the performance of the Fast Cadzow's algorithm is comparable to that of the Cadzow's algorithm but the former one is computationally more efficient.

\begin{table}[ht!]
\caption{Average computational time (seconds), average number of iterations and average mean square error of Cadzow and Fast Cadzow over 10 random problem instances per $(n,r)$ for $N\in\lbrace 4096, 128\times 128, 64\times64\times64\rbrace$ and  $r \in \lbrace 5,10,20\rbrace$.
}
\centering
\label{tab:spectral_missing_noise}
\vspace{0.2cm}
\begin{tabular}{|c|ccc|ccc|ccc|}
\hline
$r$          & \multicolumn{3}{c|}{5}    & \multicolumn{3}{c|}{10}  &\multicolumn{3}{c|}{20}  \\
\hline           & MSE    & Iter  & Time    & MSE   & Iter  &Time  & MSE    & Iter   & Time  \\
\hline           & \multicolumn{9}{c|}{$N = 4096$}                            \\
\hline Cadzow     & 4.32e-02 & 27.5 & 1.10  & 6.13e-02 & 28.5 & 1.83  & 8.40e-02 & 30.4 & 5.08\\
\hline Fast Cadzow & 4.32e-02 & 27.5 & 0.43  & 6.13e-02 & 28.5 & 0.95 & 8.40e-02 & 30.5 & 2.54 \\
\hline           & \multicolumn{9}{c|}{$N = 128\times128$}                            \\
\hline Cadzow     & 2.86e-02 & 27.0 & 1.13  & 4.02e-02 & 27.6 & 1.77 &5.91e-02 &28.8 & 4.35 \\
\hline Fast Cadzow & 2.86e-02 & 27.0 & 0.50  &4.02e-02 & 27.7 & 0.84 &5.91e-02 &29.2 & 1.84\\
\hline                      & \multicolumn{9}{c|}{$N = 64\times64\times64$}                                                         \\
\hline Cadzow               &1.01e-02 & 20.3 & 15.33  & 1.38e-02 & 20.8 & 24.55  & 1.97e-02 & 21.0 & 53.39  \\
\hline Fast Cadzow           &1.01e-02 & 20.3 & 6.95  & 1.38e-02 & 20.8 & 13.50  &1.97e-02 & 21.0 & 28.23  \\
\hline         
\end{tabular}
\end{table}
\subsection{Denoising in the reconstruction of periodic stream of Diracs }\label{subsec:numerics2}

\begin{figure}[ht!]
\centering
\includegraphics[width=0.9\textwidth]{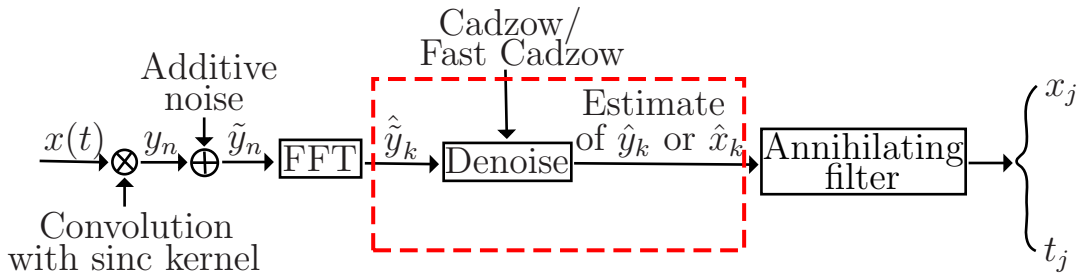}
\caption{Sampling and reconstruction of $1$-periodic stream of Diracs.}\label{fig:diracsrecon}
\end{figure}

Methods for efficient signal acquisition and reconstruction are fundamental
in signal processing. In \cite{blu2008sparse}, a novel paradigm has been proposed for sampling and reconstruction of a $\tau$-periodic stream of Diracs.  The new paradigm can achieve the minimum samples at the signal's rate of innovation. The overall sampling and reconstruction procedure is illustrated in Figure~\ref{fig:diracsrecon}.

For simplicity consider the $1$-periodic stream of Diracs $x(t)$ described in \eqref{eq:dirac_signal}.  Noticing that $x_j$ and $t_j$ ($j=1,\dots,r$) are the only $2r$ unknowns in $x(t)$, the goal  is to devise a sampling and reconstruction scheme such that $x_j$ and $t_j$ can be retrieved from about $2r$ samples. In \cite{blu2008sparse}, the signal is first convolved with a   sinc kernel of bandwidth $B$ where $B$ is an odd integer and then the samples are obtained at a equally-spaced grid. This can be formally expressed as 
\begin{align}
    y_n = \langle x(t),\mbox{sinc}(B(nT-t))\rangle = \sum\limits_{j=1}^rx_j\phi(nT-t_j),\quad  n=1,2,\cdots N,
\end{align}
where $T =1 /N$, and $\phi(t)$ is the Dirichlet kernel:
\begin{align}
    \phi(t) = \frac{\sin(\pi Bt)}{B \sin(\pi t)}.
\end{align}

A  annihilating filter method has been proposed to reconstruct $x(t)$ from the $N$ samples $\{y_n\}$. A key ingredient in the method is the construction of a filter which can annihilate $e^{-i2\pi t_j}$. Recall that $x(t)$ is a periodic signal with discretized Fourier coefficients $\{\hat{x}_k\}$. It has been shown in \cite{blu2008sparse} that the annihilation filter can be identified as a null space vector of a Toeplitz matrix involving a  set of consecutive values  of $\hat{x}_k$. Then the reconstruction of $x(t)$ from the samples finally boil down to the problems of estimating $\hat{x}_k$ from $y_n$. To this end, it is shown in \cite{blu2008sparse} that the Fourier coefficients $\{\hat{y}_k\}$ of $\{y_n\}$ coincide with a subset of $\{\hat{x}_k\}$.

If there is no noise, the annihilating filter method is able to reconstruct $x(t)$ exactly from the minimum number of measurements. However, noise is prevalent during the sampling process. We will consider the additive noise model; namely, $\tilde{y}_n=y_n+\mbox{noise}$ is observed. When noise exists a common strategy 
is to oversample the measurements and then do denoising. This is where the Cadzow's or Fast Cadzow's algorithm plays a role. In Section~\ref{subsec:intro1} we have shown that the Hankel matrix corresponding to a subset of $\hat{x}_k$ is low rank. Noticing the agreement between $\{\hat{y}_k\}$ and $\{\hat{x}_k\}$, it follows immediately that the Hankel matrix associated with 
$\{\hat{y}_k\}$ is also low rank. Therefore,  the Cadzow's or Fast Cadzow's algorithm can be used for the task of denoising.

\begin{table}[ht!]
\caption{Average computational time (seconds), average number of iterations and average mean square error of Cadzow and Fast Cadzow over 1500 random problem instances for $\epsilon \in \lbrace 0.1, 0.3, 0.5\rbrace$. }
\centering
\label{tab:dirac_denosing}
\vspace{0.2cm}
\begin{tabular}{|c|ccc|ccc|ccc|}
           \hline  $\epsilon$ & \multicolumn{3}{c|}{0.1}   & \multicolumn{3}{c|}{0.3}   & \multicolumn{3}{c|}{0.5}   \\ \hline
           & MSE    & Iter    & Time   & MSE   & Iter    & Time   & MSE    & Iter    & Time   \\ \hline
Cadzow     & 5.62e-02 & 16.35 & 0.0234 & 1.79e-01 & 17.57 & 0.0251 & 3.29e-01 & 18.43 & 0.0266 \\ \hline
Fast Cadzow & 5.62e-02 & 16.42 & 0.0055 & 1.79e-01 & 17.69 & 0.0060 & 3.29e-01 & 18.64 & 0.0065\\ \hline
\end{tabular}
\end{table}

The denoising problem studied here differs from the one in the last subsection as the Hankelization of the signal is not low rank but the Hankelization of its Fourier coefficients is low rank. Next we compare the performance of the Cadzow's and Fast Cadzow's algorithms in this situation following the setup from \cite{blu2008sparse}. Tests are conducted for  $x(t)$ with $r=7$, where $x_j$ is uniformly sampled from $[0.5,1.5]$ and $t_j$ is uniformly sampled from $[0,1)$. The window length of the sinc kernel is chosen to $B=71$ and $N=B$ samples are observed.
We test three different noise levels and 1500 randomly generated
problem instances are simulated for each noise level. The algorithms are terminated when $\|\bz_{k+1}-\bz_k\|/\|\bz_k\|\leq 10^{-6}$. The average computational results are shown in Table~\ref{tab:dirac_denosing}. Even though both of the algorithms are very fast in the experiments due to the small problem size, the Fast Cadzow's algorithms is about four times faster. 
\subsection{Seismic data denoising and reconstruction}\label{subsec:numerics3}
Seismic  denoising and seismic recovery from missing traces are two major tasks in seismic data processing and different techniques have been developed. The Cadzow's algorithm which was called the $f$-$x$ eigen filtering has been proposed to  attenuate random noise from seismic records by Trickett \cite{sacchi2009fx,oropeza2009multifrequency,trickett2008f,trickett2009prestack,oropeza2011simultaneous}. It is based on the idea that if seismic data consists of $r$ linear events  then the associated Hankel matrices of the constant-frequency slices are of  rank $r$. The algorithm was extended in \cite{oropeza2011simultaneous} to the problem of reconstructing missing traces (see \eqref{eq:cadzow_missing}), which can be viewed as a non-convex version of   the method of projection
onto convex sets (POCS, \cite{abma20063d}). Moreover, a generalized version of the algorithm was also presented in \cite{oropeza2011simultaneous} for  simultaneously denoising and reconstruction  of the seismic data, which is actually  the one listed in \eqref{eq:cadzow_denoising_missing}.

\begin{figure}[ht!]
    \centering
    \includegraphics[width = 4cm,height = 8cm]{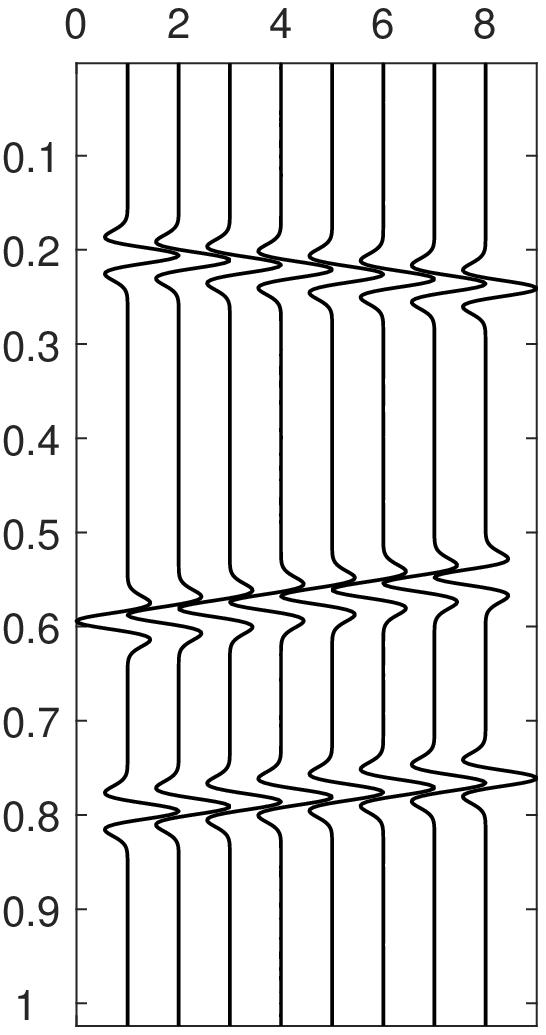}
    \includegraphics[width = 4cm,height = 8cm]{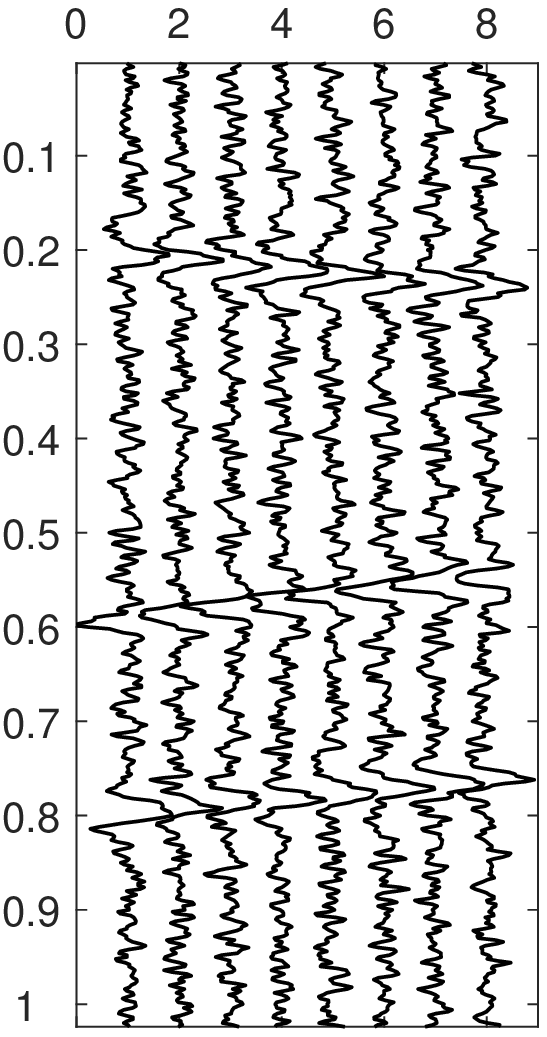}
    \includegraphics[width = 4cm,height = 8cm]{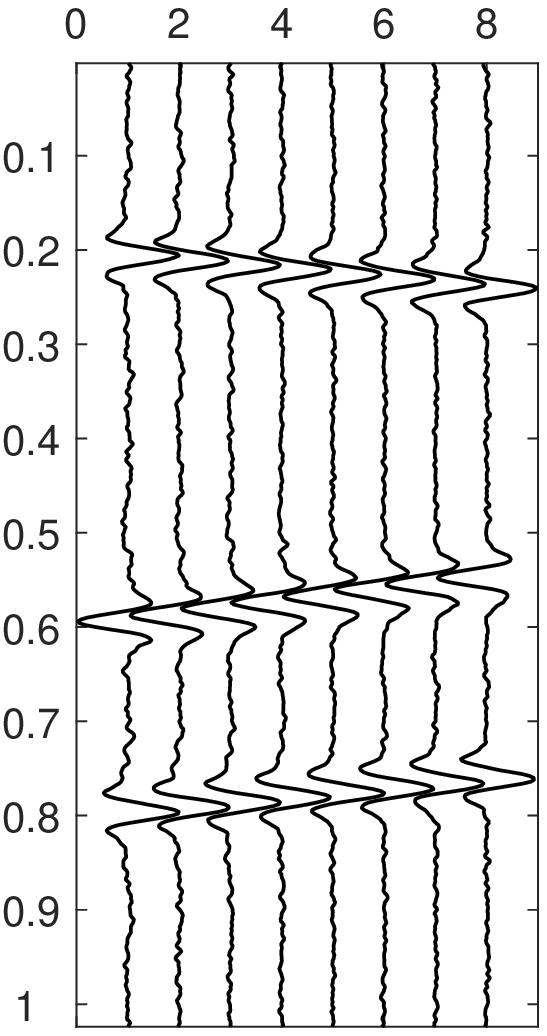}
    \includegraphics[width = 4cm,height = 8cm]{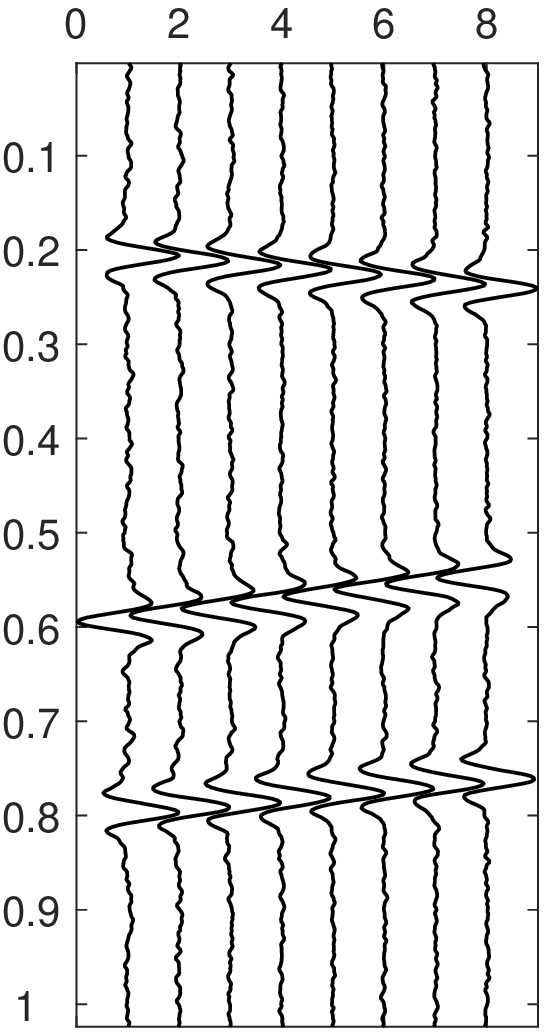}
    \caption{A slice of (a) noiseless data, (b) noisy data, (c) data returned by Cadzow denoising, and (d) data returned Fast Cadzow denoising. Vertical: time dimension; Horizontal: one spatial dimension with the other three being held fixed.}
    \label{fig:my_labelx}
\end{figure}

\begin{figure}[ht!]
    \centering
    \includegraphics[width = 4cm,height = 8cm]{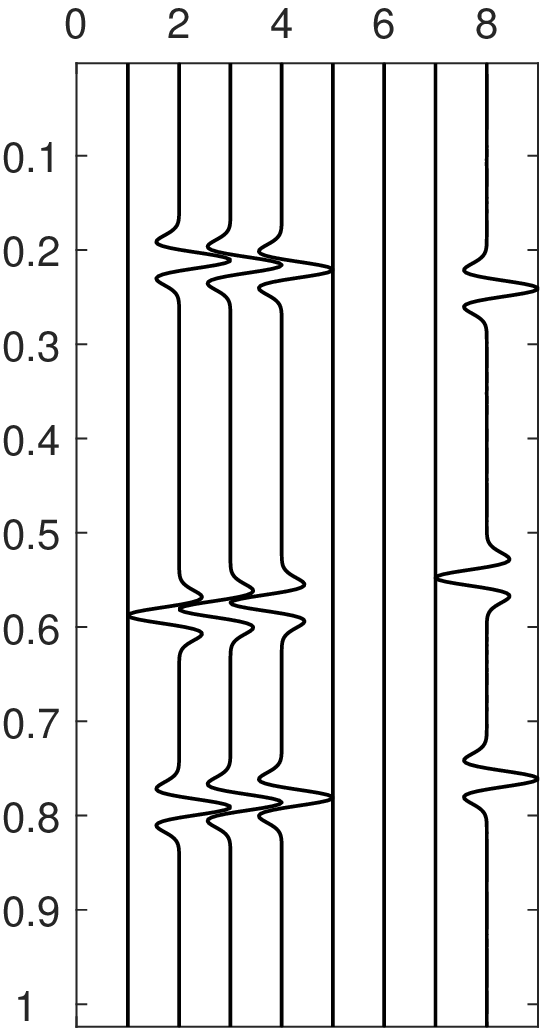}
    \includegraphics[width = 4cm,height = 8cm]{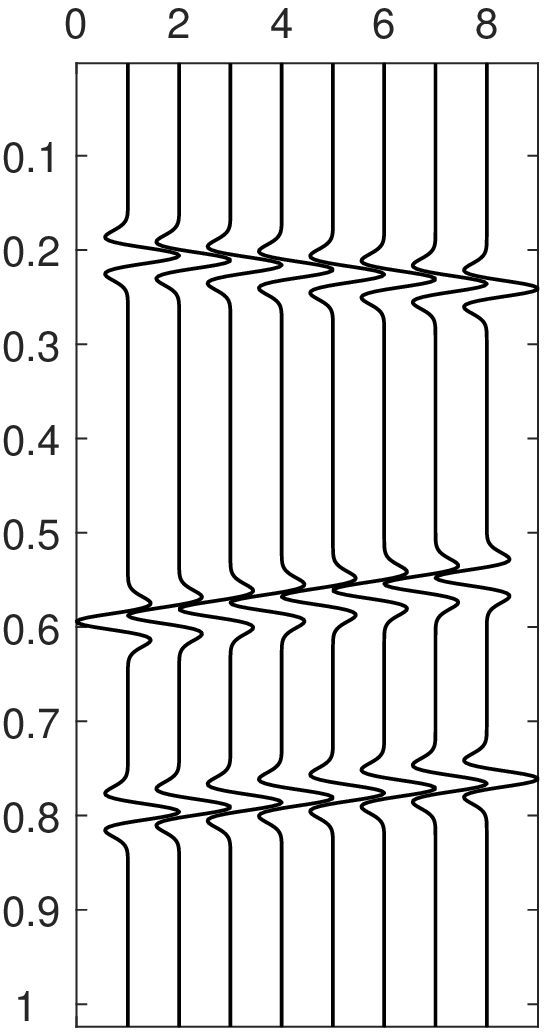}
    \includegraphics[width = 4cm,height = 8cm]{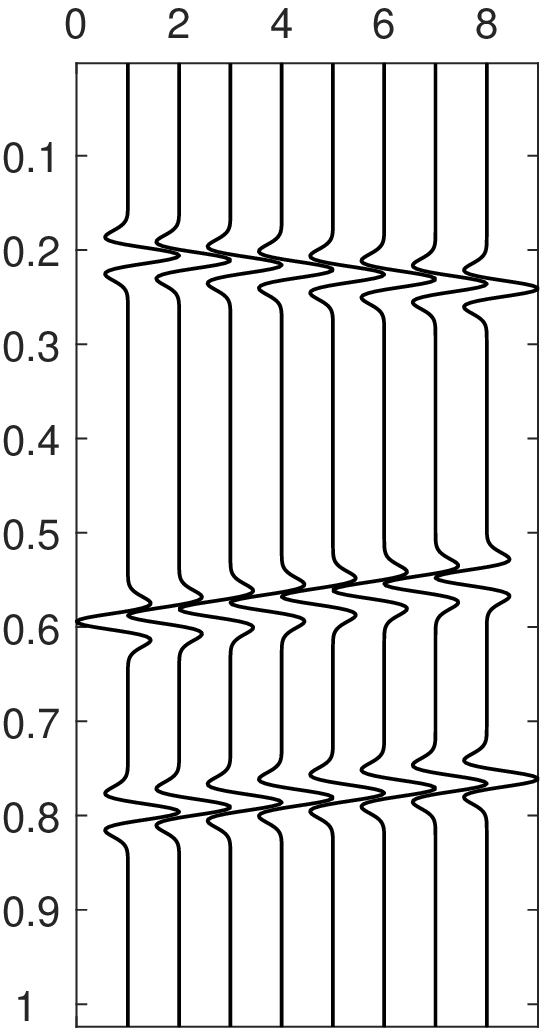}
    \caption{A slice of (a) data with missing traces, (b) data returned by Cadzow recovery, and (c) data returned by Fast Cadzow recovery. Vertical: time dimension; Horizontal: one spatial dimension with the other three being held fixed.}
    \label{fig:my_labely}
\end{figure}

As an accelerated variant of the Cadzow's algorithm, the Fast Cadzow's algorithm is equally applicable for seismic data processing.  We compare the performance of the two algorithms on a 4D volume with $8\times 8\times 8\times 8$ traces and each trace contains $512$ time samples. Thus, this is a  5D data array with one time dimension and four spatial dimensions.
The data consists of three linear events, so we set $r=3$ in the algorithms.
Tests have been conducted for the denoising problem with $\epsilon=1$ as well as for the recovery problem with 50\% missing traces. As is typical in the literature on seismic data processing, the algorithms are run for a total number of $10$ iterations. The MSE and computational time are presented in Table~\ref{tab:my-table6}. In addition, the seismic profile of one data slice with three fixed spatial dimensions has been plotted in Figure~\ref{fig:my_labelx} and \ref{fig:my_labely} for seismic denoising and recovery, respectively. From the table and the figures we can see that the two algorithms have comparable performance. However, the table does show that the MSE of the Fast Cadzow's algorithm is slightly better than that of  the Cadzow's algorithm for seismic denoising. Moreover, the Fast Cadzow's algorithm is two times faster for the denoising problem and four times faster for the recovery problem.

\begin{table}[ht!]
\caption{Mean square error and computational time (seconds) of Cadzow and Fast Cadzow over seismic denoising and seismic recovery from missing traces.  }
\vspace{.2cm}
\centering
\label{tab:my-table6}
\begin{tabular}{|c|cc|cc|}
\hline           & \multicolumn{2}{c|}{Seismic denoising} & \multicolumn{2}{c|}{Seismic recovery} \\
\hline           & MSE         & Time        & MSE         & Time          \\
\hline Cadzow     & 9.58e-02      & 78.34    & 9.88e-04      & 139.13      \\
\hline Fast Cadzow & 8.88e-02      & 37.08     & 9.88e-04      & 34.76    \\ \hline   
\end{tabular}
\end{table}

\section{A gradient variant for denoising}\label{sec:gradient}
\subsection{Motivation and new algorithms}
\begin{figure}[ht!]
    \centering
    \includegraphics[width=0.45\textwidth]{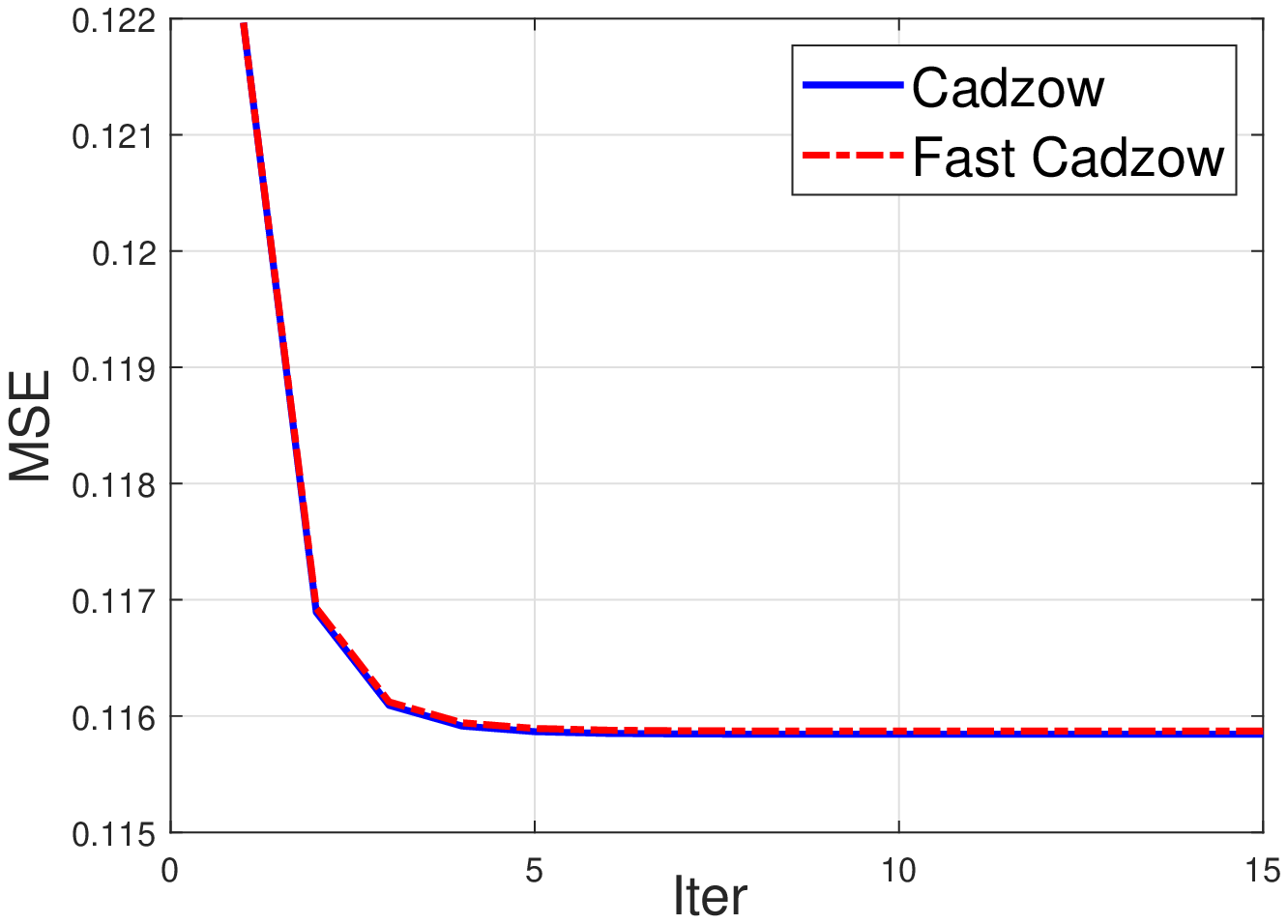}
    \includegraphics[width =0.45\textwidth]{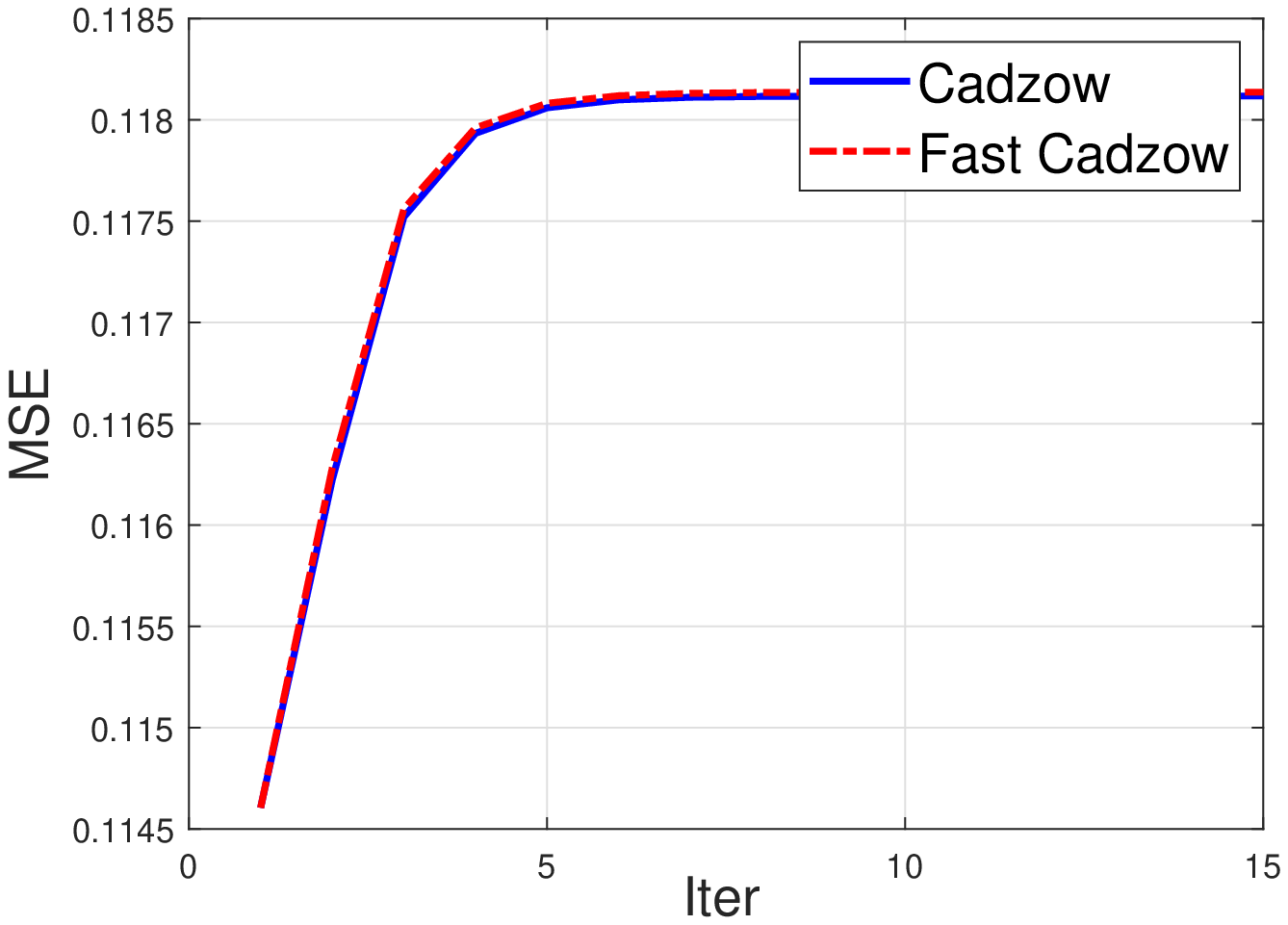}
    \caption{ (a) MSE decreasing  case, (b) MSE increasing case.}\label{fig:grad_fg1}
\end{figure}

When applying the Cadzow's algorithm or the Fast Cadzow's algorithm for signal denoising, it has been observed that the MSE does not always decrease as the algorithm iterates. For example, in \cite{trickett2008f} Trickett notes that 
\begin{quote}\em
Cadzow recommended iterating between the rank reduction and averaging steps, but I have not found it to give better results for this application.
\end{quote}
Gillard also notes in his paper \cite{gillard2010cadzow} that
\begin{quote}\em
In the simulation study within this paper, it has been
demonstrated that repeated iterations of Cadzow’s basic algorithm (in an attempt to separate the noise from the signal)
may result in an increased RMSE from the true signal. 
\end{quote}
As an illustration, the MSE plots of the Cadzow's and Fast Cadzow's algorithms corresponding to two random instances on 1D spectrally sparse signal denoising with $N=256$, $r=5$, and $\epsilon=0.5$ are presented in Figure~\ref{fig:grad_fg1}. For both algorithms, we can see that the MSE in the left plot  decreases but the MSE in the right plot increases. In other words, it is possible that the MSE of the two algorithms may increase for signal denoising, which motivates us to find some new alternatives.

\begin{figure}[ht!]
\centering
\includegraphics[width=0.45\textwidth]{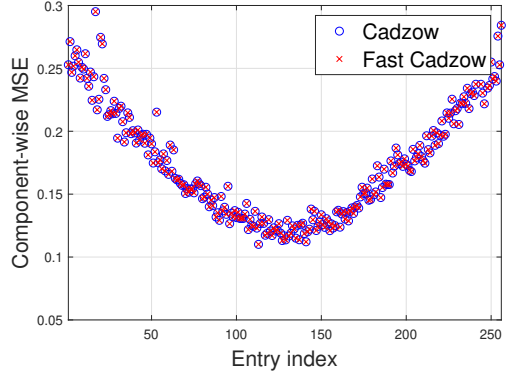}
\caption{The average component-wise MSE  after the Cadzow's and Fast Cadzow's denoising.}\label{fig:grad_fg2}
\end{figure}

To explain this phenomenon, we will study the equivalent form of the Cadzow's algorithm in \eqref{eq:cadzow_matrix}. It is trivial that the update can be rewritten as 
\begin{align*}
\BZ_{k+1} &= \P_{\M_H}\P_{\M_r}\BZ_k\\
&=\P_{\M_H}\P_{\M_r}(\BZ_k+t(\BY-\BZ_k)),\quad t=0,
\end{align*}
where $\BZ_k=\H\bz_k$ and $\BY = \H\by=\H(\bx+\be)$. Thus, the Cadzow's algorithm can be interpreted as a projected gradient method for the following optimization problem
\begin{align*}
\min\frac{1}{2}\|\BZ-\BY\|_F^2\quad \mbox{subject to}\quad \rank(\BZ)\leq r \mbox{ and }\BZ\mbox{ is Hankel},\numberthis\label{eq:cadzow_opt}
\end{align*}
though with a step length $t=0$. Note that the objective function in the above optimization problem is equal to 
\begin{align*}
\sum_{a=0}^{N-1}w_a|z_a-y_a|^2,
\end{align*}
where $w_a$ is the number of entries on the $a$-th skew-diagonal
of an $L\times K$ matrix; see \eqref{eq:weight}. Thus, the Cadzow's algorithm for signal denoising indeed solves an optimization problem which puts different weights on different entries of the signal. Moreover, the Fast Cadzow's algorithm should share the same property as the Cadzow's algorithm.
Since the weights for the middle entries are larger than those at the two ends, one may expect that after the Cadzow's or Fast Cadzow's denoising the component-wise MSE (defined as $|z_k(i)-x(i)|/|x(i)|,~i=0,\cdots,N-1$) of the middle entries should be smaller.
The plot in Figure~\ref{fig:grad_fg2} confirms this speculation, which shows the average component-wise MSE out of 1500 random tests after applying the Cadzow's and Fast Cadzow's algorithms for 1D spectrally sparse signal denoising with $N=256$, $r=5$, and $\epsilon=0.5$.

In order to address the unbalanced weight issue, we consider the following optimization problem with re-weighted objective function:
\begin{align*}
\min\frac{1}{2}\|\sqrt{\BW}\odot(\BZ-\BY)\|_F^2\quad \mbox{subject to}\quad\rank(\BZ)\leq r \mbox{ and }\BZ\mbox{ is Hankel},\numberthis\label{eq:gradient_opt}
\end{align*}
where $\sqrt{\BW}$ means taking the square-root of each entry of $\BW$,  $\odot$ denotes the component-wise product and 
\begin{align*}
\sqrt{\BW}=\H(\sqrt{w_0},\cdots,\sqrt{w_{N-1}})= \begin{bmatrix}
1 & \frac{1}{\sqrt{2}} & \frac{1}{\sqrt{3}} & \vdots & \vdots\\
\frac{1}{\sqrt{2}} & \frac{1}{\sqrt{3}}& \vdots & \vdots & \vdots\\
\frac{1}{\sqrt{3}} & \vdots & \vdots & \vdots &\frac{1}{\sqrt{3}}\\
\vdots & \vdots & \vdots & \frac{1}{\sqrt{3}} & \frac{1}{\sqrt{2}}\\
\vdots & \vdots & \frac{1}{\sqrt{3}} &\frac{1}{\sqrt{2}}& 1
\end{bmatrix}.
\end{align*}

\begin{algorithm}[ht!]
\caption{Gradient Method}
\label{alg:grad}
\begin{algorithmic}[1]
\Require $\bz_0=\by$
\For{$k=0,1,\cdots$}
\State $\bz_{k+1}=\H^\dagger\T_r\H\lb\bz_k+\frac{1}{\BD}\odot(\by-\bz_k)\rb$
\EndFor
\Ensure $\bz_k$
\end{algorithmic}
\end{algorithm}

A projected gradient method with step length $t=1$ can then be developed for \eqref{eq:gradient_opt} as follows: 
\begin{align*}
\BZ_{k+1} = \P_{\M_H}\P_{\M_r}(\BZ_k+\BW\odot(\BY-\BZ_k)).
\end{align*}
The equivalent form of the algorithm in the vector domain, dubbed Gradient method, is presented in Algorithm~\ref{alg:grad}, where $1/\BD = [1/w_0,\cdots,1/w_{N-1}]^\top$.

\begin{algorithm}[ht!]
\caption{Fast Gradient Method}
\label{alg:fastgrad}
\begin{algorithmic}[1]
\Require $\bz_0=\by$
\For{$k=0,1,\cdots$}
\State Compute a matrix subspace $T_k$
\State $\bz_{k+1}=\H^\dagger\T_r\P_{T_k}\H\lb\bz_k+\frac{1}{\BD}\odot(\by-\bz_k)\rb$
\EndFor
\Ensure $\bz_k$
\end{algorithmic}
\end{algorithm}

Of course we can use the same technique as in the Fast Cadzow's algorithm to accelerate Algorithm~\ref{alg:grad}, leading to the Fast Gradient method; see Algorithm~\ref{alg:fastgrad}. In the remainder of this section we compare 
the  algorithms on different signal denoising problems, focusing on the denoising performance of the algorithms. About the computation efficiency,
it is clear that the dominant per iteration computational cost of the Gradient method is the same as that of the Cadzow's algorithm, and the dominant per iteration computational cost of the Fast Gradient method is the same as that of the Fast Cadzow's algorithm.
\subsection{Empirical evaluations}
\begin{figure}[ht!]
    \centering
    \includegraphics[width=0.3\textwidth]{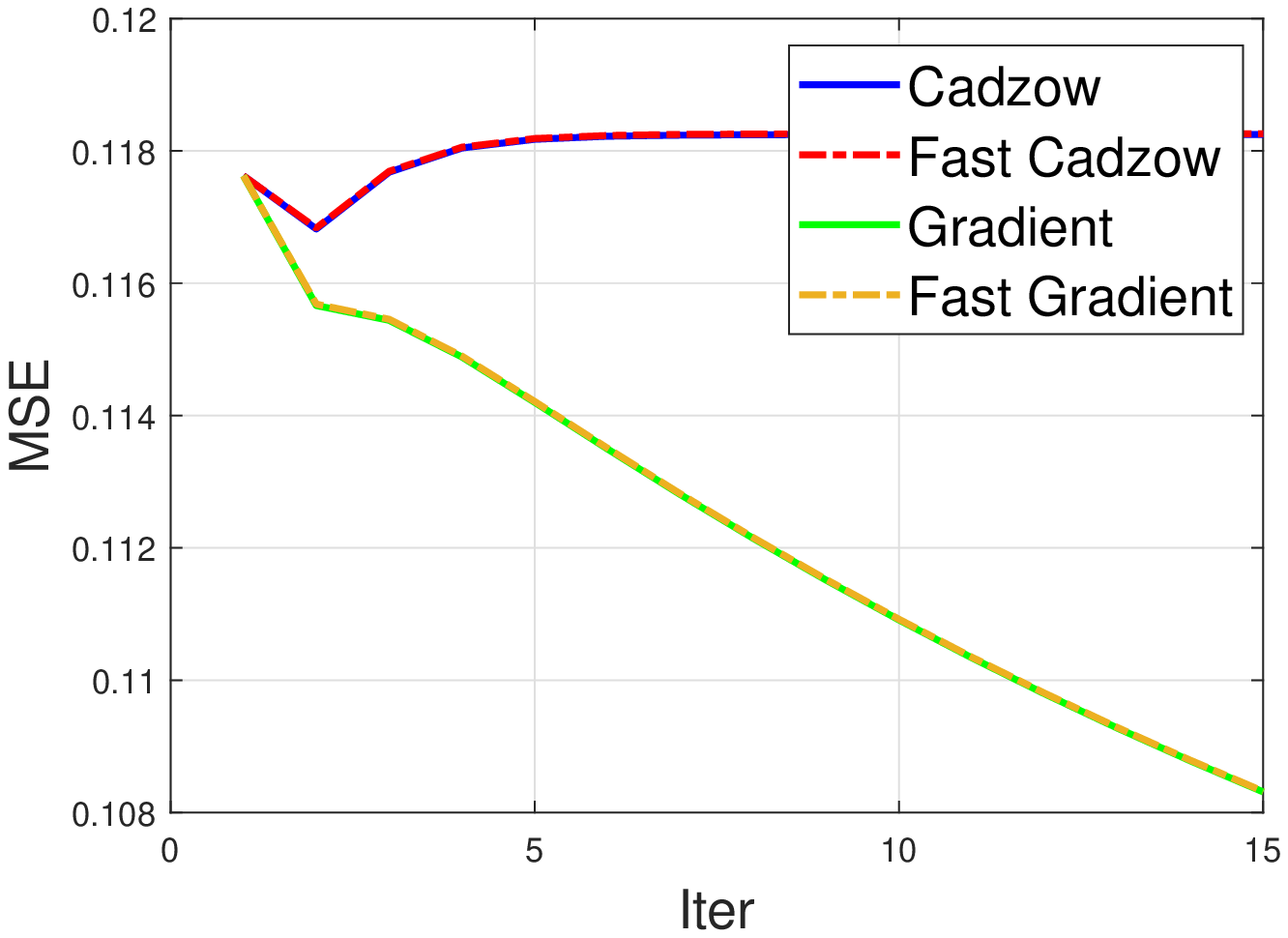}
    \includegraphics[width=0.3\textwidth]{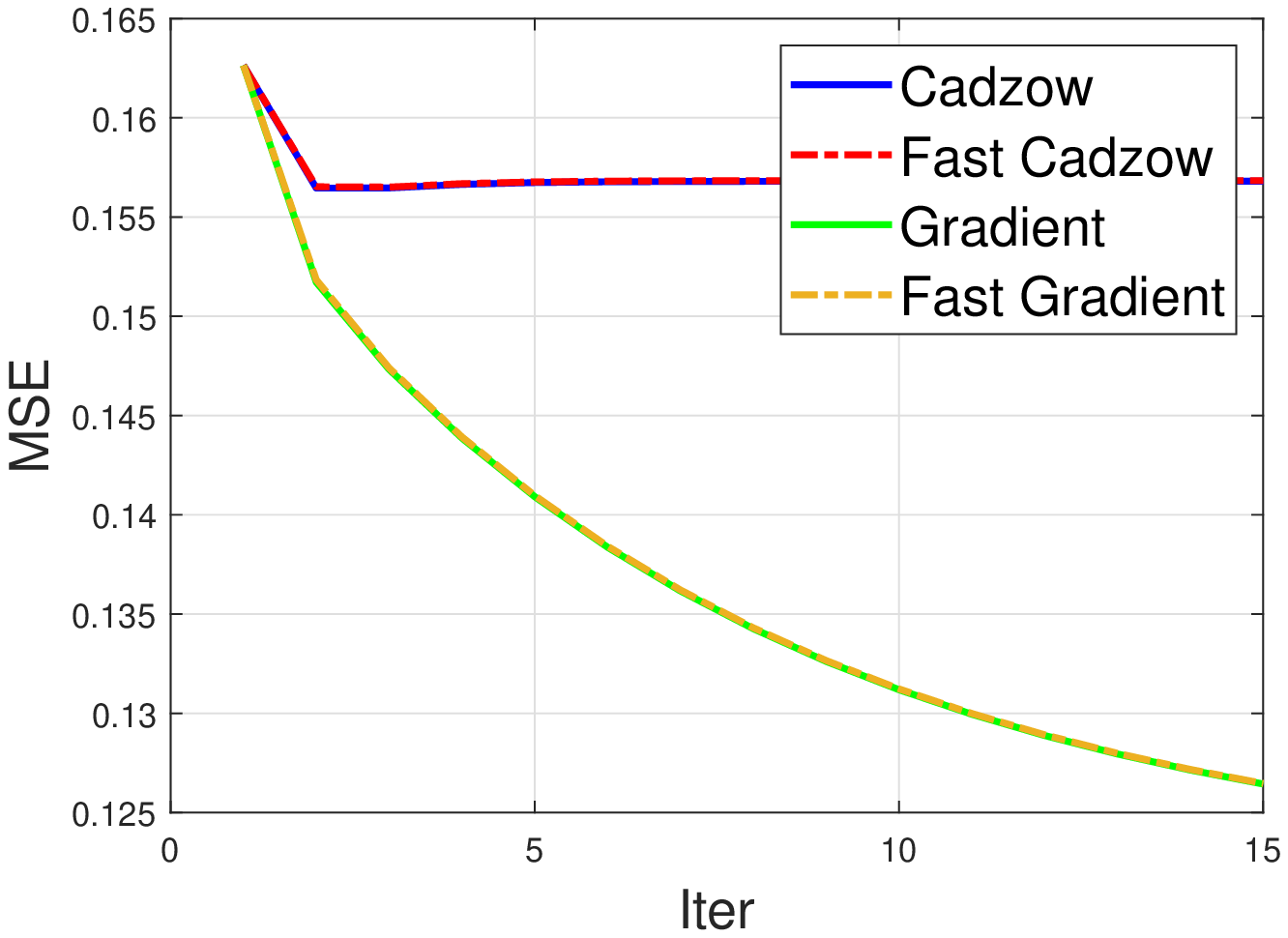}
    \includegraphics[width=0.3\textwidth]{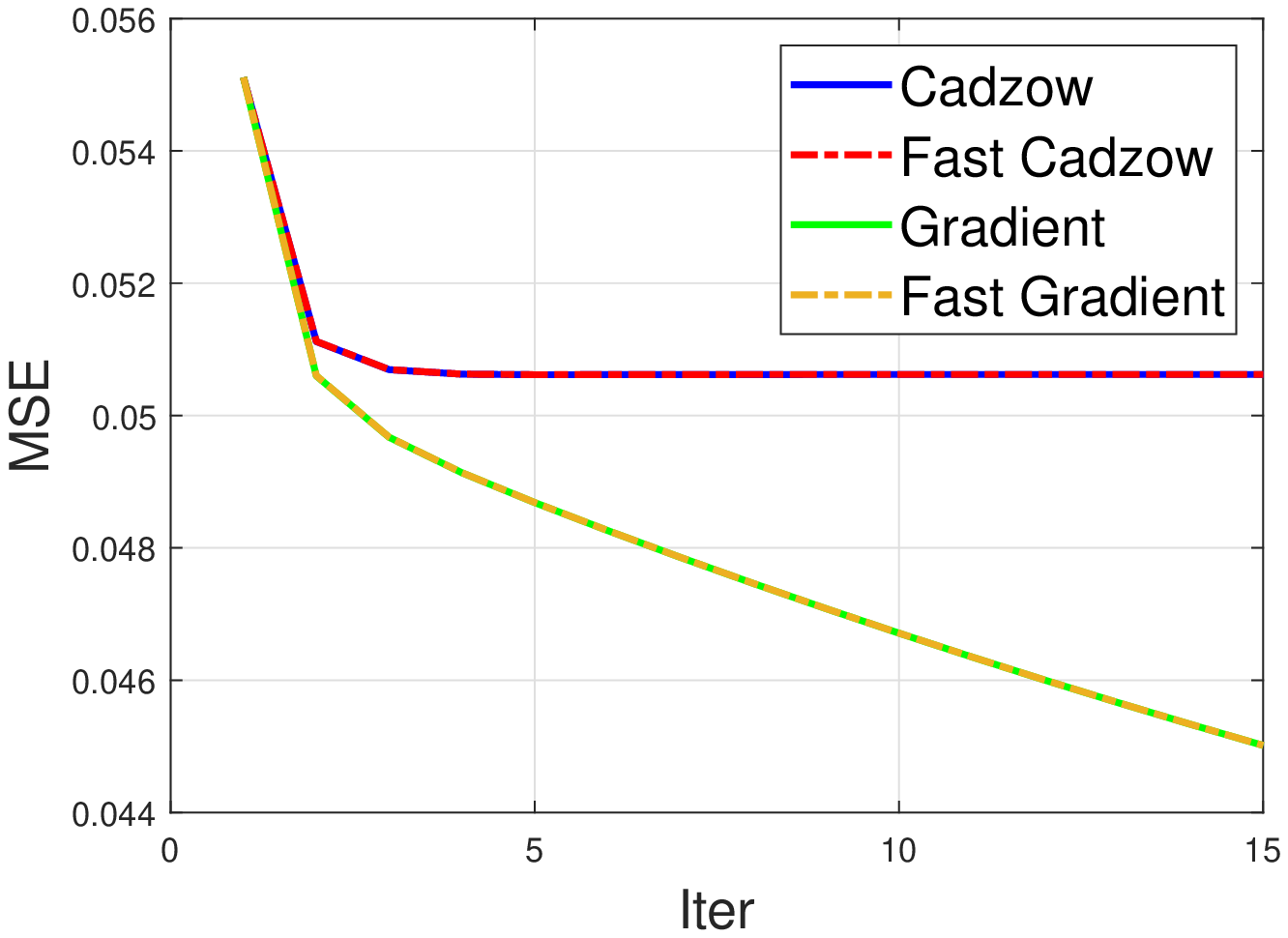}
    \caption{Average per iteration MSE out of $1500$ random tests for spectrally sparse signal denoising.}
    \label{fig:grad_fg3}
\end{figure}

We first compare the four algorithms, Cadzow, Fast Cadzow, Gradient and Fast Gradient, on spectrally sparse signal denoising (see Section~\ref{subsec:numerics1}) with $N\in\lbrace 256, 16\times 16, 16\times16\times16\rbrace$, $r = 5$ and noise level $\epsilon = 0.5$.
All the algorithms are run for the maximum of $15$ iterations. For each $N$, 1500 random problem instances are tested. The average per iteration MSE of each algorithm against the number of iterations is presented in Figure~\ref{fig:grad_fg3}. Overall the MSE of the Cadzow's and Fast Cadzow's algorithms for the 1D denoising problem increases as the algorithm iterate. In contrast, the MSE of the Gradient and Fast Gradient methods decreases as the algorithms iterate. For the 2D and 3D denoising problems, though on average the MSE of all the algorithms shows a decreasing trend, the Gradient and Fast Gradient methods can achieve smaller MSE.

In addition, we also examine each individual test and consider one to be a positive test if the MSE of an algorithm at the last iteration is smaller than that at the first iteration. The portion of positive tests for each test algorithm is presented in Table~\ref{tab:grad_tb1}. It can be observed that the number of  positive tests of the Gradient and Fast Gradient methods is larger than that of the Cadzow's and Fast Cadzow's algorithms, especially for the 1D denoising problem.

\begin{table}[ht!]
\centering
\caption{Portion of positive tests for the four test algorithms.}
\label{tab:grad_tb1}
\vspace{0.2cm}
\begin{tabular}{|c|c|c|c|}
 \hline             & N = 256 & N =16$\times$16 & N = 16$\times$16$\times$16 \\ \hline
Cadzow        & 0.4200  & 0.8300   & 0.9907       \\ \hline
Fast Cadzow   & 0.4193  & 0.8293   & 0.9907       \\ \hline
Gradient      & 0.9947  & 1.0000   & 1.0000       \\ \hline
Fast Gradient & 0.9947  & 1.0000   & 1.0000 \\ \hline    
\end{tabular}
\end{table}

\begin{figure}[ht!]
\centering
\includegraphics[width=0.45\textwidth]{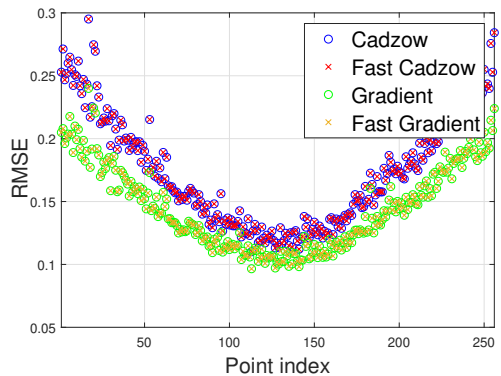}
\caption{The average component-wise MSE after  denoising.}
\label{fig:grad_fg4}
\end{figure}

It is also interesting to compare the average component-wise MSE of the outputs of  the four algorithms; see Figure~\ref{fig:grad_fg4} for the 1D case. The figure shows that, to some extent, the Gradient and Fast Gradient methods are able to alleviate the issue of unbalanced weights over different entries.

\begin{figure}[ht!]
\centering
\includegraphics[width=0.45\textwidth]{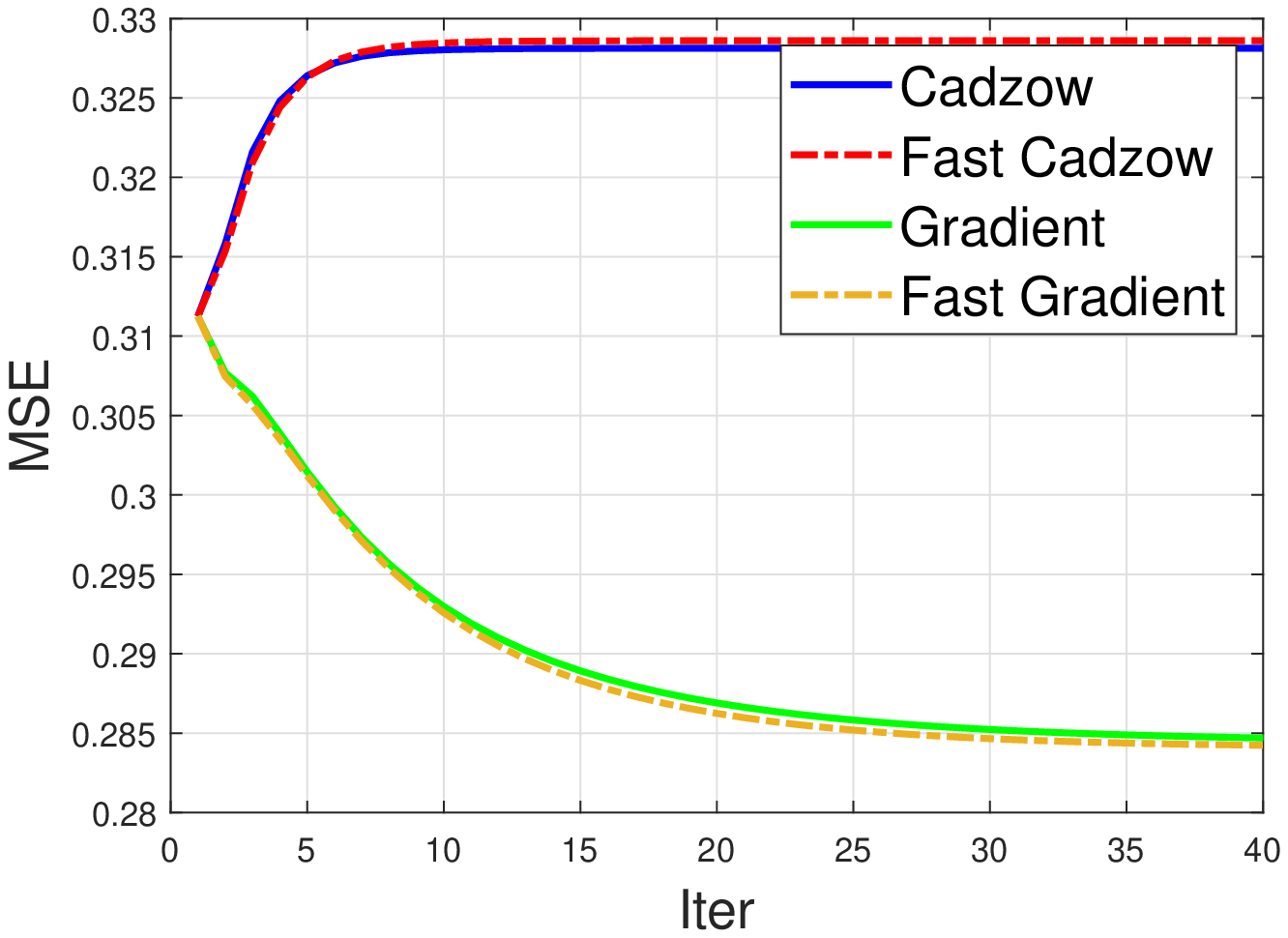}
\includegraphics[width=0.45\textwidth]{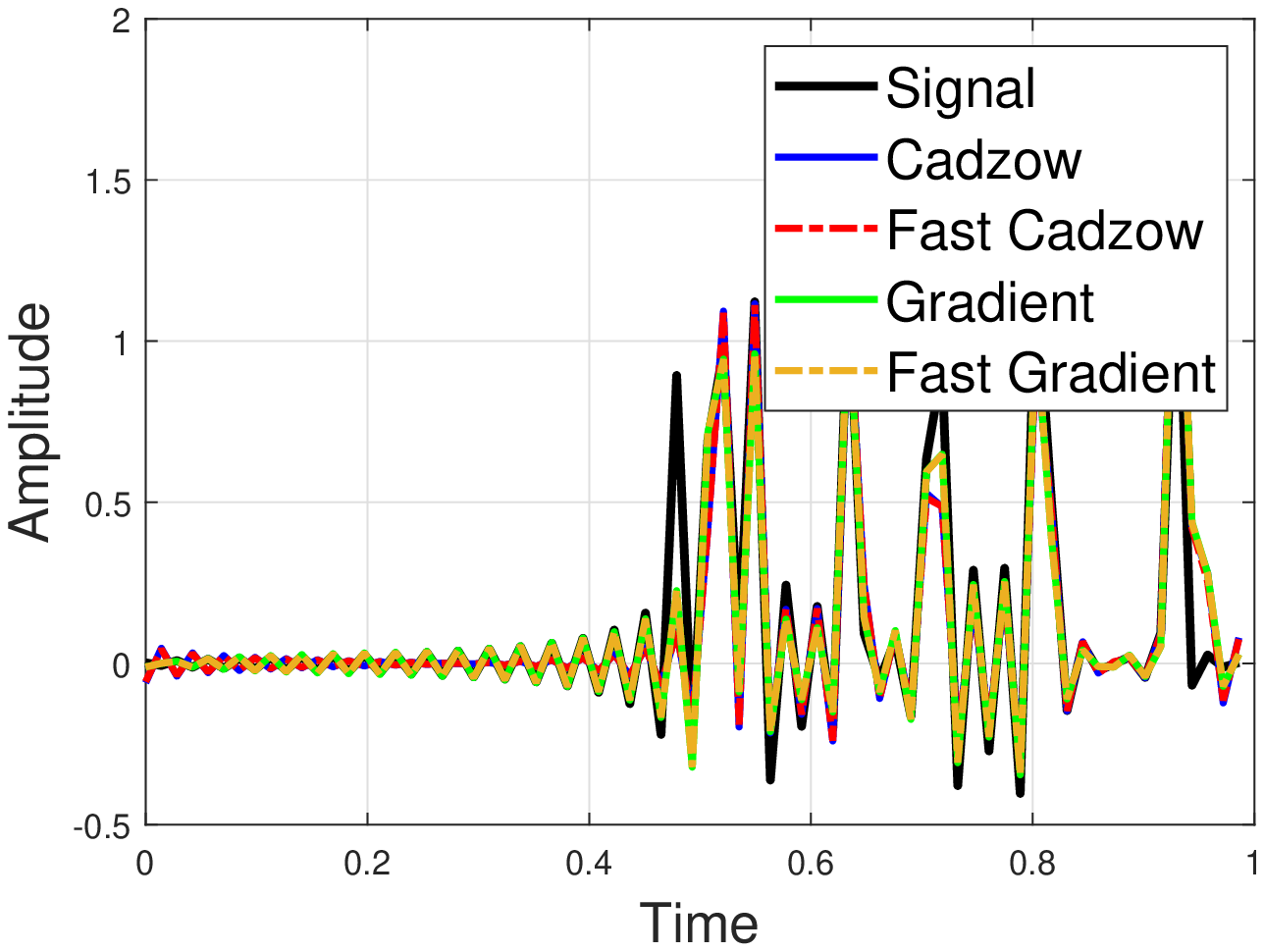}
\caption{Left: Average per iteration MSE out of $1500$ random tests for denoising Dirac samples in the Frequency domain; Right: Denoising results in the time domain.}\label{fig:grad_fg5}
\end{figure}

Next we test the algorithms for the denoising problem arising from the reconstruction of Diracs; see Section~\ref{subsec:numerics2}. Roughly speaking, the Fourier transform of the test signal has a low rank Hankelization. Thus we first do denoising  in the Frequency domain and then get the denoised signal via the inverse FFT. Here we  repeat 1500 random tests  for the noise level $\epsilon=0.5$. The average MSE against the number of iterations is plotted in Figure~\ref{fig:grad_fg5}. Meanwhile, the portion of positive tests of each algorithm is contained in Table~\ref{tab:grad_tb2}. Obviously, the average MSE of the Gradient and Fast Gradient methods decreases while the MSE of the Cadzow's and Fast Cadzow's algorithms increases, and there are more positive tests after the Gradient and Fast Gradient denoising. Moreover, Figure~\ref{fig:grad_fg5} also shows the denoising results of the four different algorithms in the time domain from a single random test. It can be seen that the Gradient and Fast Gradient methods exhibit better denoising performance in the region
where  $t<0.4$.

\begin{table}[ht!]
\centering
\caption{Portion of positive tests for the four test algorithms.}
\label{tab:grad_tb2}
\vspace{0.2cm}
\begin{tabular}{|c|c|c|c|c|}
 \hline      & Cadzow & Fast Cadzow & Gradient & Fast Gradient \\ \hline
N = 71 & 0.3013 & 0.2887      & 0.8167   & 0.8240   \\ \hline
\end{tabular}
\end{table}

\begin{figure}[ht!]
\centering
\includegraphics[width=0.45\textwidth]{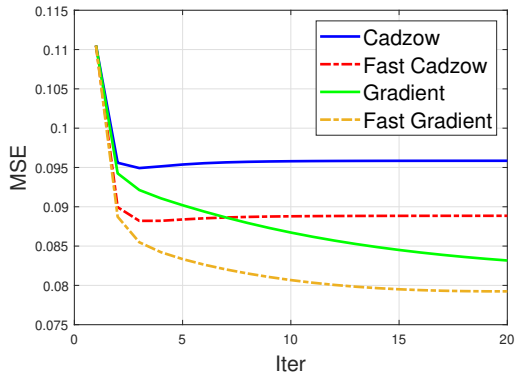}
\caption{MSE for seismic data denoising.}\label{fig:grad_fg7}
\end{figure}

Lastly, we compare the algorithms on the seismic denoising problem from Section~\ref{subsec:numerics3}. The MSE plot  against the number of iterations is presented in Figure~\ref{fig:grad_fg7}.  The plot shows that overall all the algorithms have decreasing MSE, but the MSE of the Gradient and Fast Gradient methods is smaller. It is also worth noting that for seismic data denoising the accelerated algorithms have noticeable smaller MSE than their non-accelerated counterparts. 

\section{Conclusion}\label{sec:conclusion}
In this paper we consider the signal denoising and recovery problems for  signals corresponding to low rank Hankel matrices. New algorithms have been proposed which can complete the tasks more efficiently and effectively. For future work, we would like to study the  theoretical convergence analysis of the proposed algorithms under certain random models. It may also be interesting to see whether it is possible to design better adaptive re-weighting  schemes for the gradient methods.
\section*{Acknowledgments}
HW and KW would like to thank Jianjun Gao for sending us the 
5D data for the seismic denoising and recovery tests.
\bibliographystyle{siam}
\bibliography{ref}
\end{document}